\documentstyle[amscd,amssymb,verbatim,epsf]{amsart}

\begin{document}
\theoremstyle{plain}
\newtheorem{Thm}{Theorem}
\newtheorem{Cor}{Corollary}
\newtheorem{Ex}{Example}
\newtheorem{Con}{Conjecture}
\newtheorem{Main}{Main Theorem}
\newtheorem{Lem}{Lemma}
\newtheorem{Prop}{Proposition}

\theoremstyle{definition}
\newtheorem{Def}{Definition}
\newtheorem{Note}{Note}

\theoremstyle{remark}
\newtheorem{notation}{Notation}
\renewcommand{\thenotation}{}

\errorcontextlines=0
\numberwithin{equation}{section}
\renewcommand{\rm}{\normalshape}%

\title[Casimir Effect]%
   {The Casimir Effect Between Non-parallel Plates \\
  by Geometric Optics}
\author{Brendan Guilfoyle}
\address{Brendan Guilfoyle\\
          Department of Mathematics and Computing \\
          Institute of Technology, Tralee \\
          Clash \\
          Tralee  \\
          Co. Kerry \\
          Ireland}
\email{brendan.guilfoyle@@ittralee.ie}
\author{Wilhelm Klingenberg}
\address{Wilhelm Klingenberg\\
 Department of Mathematical Sciences\\
 University of Durham\\
 Durham DH1 3LE\\
 United Kingdom}
\email{wilhelm.klingenberg@@durham.ac.uk }

\author{Siddhartha Sen}
\address{Siddhartha Sen\\
 School of Mathematics\\
 University of Dublin\\
 Dublin 2\\
 Ireland\\
and\\
Department of Theoretical Physics\\
IACS\\ 
Kolkata\\
 India}
\email{sen@@maths.tcd.ie}

\date{April 24th, 2004}

\begin{abstract}
Recent work by Jaffe and Scardicchio has expressed the optical 
approximation to the Casimir effect as a sum over geometric quantities. 
The first two authors have developed a technique which uses the
complex geometry of the space of oriented affine lines in
${\Bbb{R}}^3$ to describe reflection of rays off a surface. This
allows the quantities in the optical approximation to the Casimir
effect to be calculated. 

To illustrate this we determine explicitly and in closed form the
geometric optics approximation of the Casimir force between two
non-parallel plates. By making one of the 
plates finite, we regularise the divergence that is
caused by the intersection of the planes. In the parallel plate limit
we prove that our expression reduces to Casimir's original result.
\end{abstract}

\maketitle

\section{Introduction}
Since its discovery over 50 years ago \cite{cas} the Casimir effect
has been the subject of much research. Recently the effect has been
measured in different situations \cite{bmm} and, as a result, there is
renewed interest in calculating the Casimir effect between different
shaped boundaries \cite{lamb}. A particularly elegant approximation
based on 
geometric optics has been suggested by Jaffe and Scardicchio
\cite{jas} and they have shown that between a sphere
and a plane this approximation appears to work well. 

Their approach considers the sum over all closed paths with $n$
reflections between the boundaries.  The Dirichlet Casimir energy in this
approximation is 
\begin{equation}\label{e:casapprox}
{\cal{E}}=-\frac{\hbar c}{2\pi^2}\sum_n(-1)^n\iiint\limits_{{\cal{D}}_n}\frac{\sqrt{\Delta_n(x,x)}}{l_n^3(x)}d^3x,
\end{equation}
where $l_n$ is the length and $\Delta_n$ the Van Vleck determinant of 
a closed $n$-bounce path.

In this paper we present a complex geometry scheme
\cite{gak4} which, in principle, can be used to evaluate the Casimir
effect for a range of boundaries in the approximation scheme of Jaffe
{\it et al}. The basic ingredient of our method is a description of
reflection in ${\Bbb{R}}^3$ via the space of oriented  
affine lines, which we identify with $\mbox{T}{\Bbb{P}}^1$. In this
representation, reflection consists of the combined actions of 
$\overline{\mbox{PSL}(2,{\Bbb{C}})}$ and fibre mappings on
$\mbox{T}{\Bbb{P}}^1$.  

We compute the geometric approximation to the Casimir energy by iterating the 
above representation. The convergence of the integral is then studied
through the analytic properties of the maps involved. Of particular
interest is the case where we have translational symmetry. This
reduces the computation to a planar problem which can be addressed 
using trigonometric functions.

We illustrate the method by evaluating the case of non-parallel
plates. For plane boundaries, the Van Vleck determinant is the
reciprocal of the length squared of the path, and our approach gives
an explicit iteration of the classical method of images using
properties of the above mentioned action on $\mbox{T}{\Bbb{P}}^1$.  

We first compute the integrand of (\ref{e:casapprox}) and reduce the
infinite sum over closed paths to a finite sum. Let ($R$, $\psi$) be polar coordinates in the plane, where $\psi$ is measured from the vertical. Then: 

\vspace{0.2in}
\noindent {\bf Main Theorem 1}.

{\it 
There exists a closed $2m$-bounce path in a wedge of opening angle $\gamma$
iff $\gamma<\frac{\pi}{2m}$. There exists at most two closed paths 
(traversed in opposite directions) and the total length of a closed $2m$-bounce
reflection path  starting at the point ($R\sin\psi$, $R\cos\psi$) is
$2R|\sin(m\gamma)|$. 
 
There exists a closed $2m+1$-bounce path starting at the
point ($R\sin\psi$, $R\cos\psi$) iff $\gamma<
\frac{\pi}{2m}$ and either of the following hold:
$\psi<\pi-(m+1)\gamma$ or $\psi>m\gamma$.

There exists exactly none, one or two closed $2m+1$-bounce paths from a given point,
according to the inequalities above.

The total length of a closed $2m+1$-bounce path is
$2R|\cos[\psi+(m+1)\gamma]|$ in the former case, and
$2R|\cos(\psi-m\gamma)|$ in the latter case.  
} 
\vspace{0.2in}

There remains a divergence in the expression for the Casimir energy 
caused by the vertex of the wedge. We remove this by restricting one
of the plates to be finite.  This introduces
restrictions on the domain of integration ${\cal{D}}_n$ which must be
carefully dealt with separately in the even and odd bounce cases. 

Consider a plate lying between radial
distances $R=R_0$ and $R=R_1$ from the origin, forming an angle
$\gamma$ with the horizontal plane containing the origin. Define $m_0$
and $m_1$ by
$\frac{\pi}{2m_0+2}\le\gamma<\frac{\pi}{2m_0}$, and either
$\cos(m_1\gamma)\le\frac{R_0}{R_1}\le
\frac{\cos[(m_1-1)\gamma]}{\cos\gamma}$ for $m_1$ even, or
$\frac{\cos(m_1\gamma)}{\cos\gamma}\le\frac{R_0}{R_1}\le\cos[(m_1-1)\gamma]$ 
for $m_1$ odd. Our result is:

\vspace{0.2in}

\noindent {\bf Main Theorem 2}.

{\it 
The geometric optics approximation to the Casimir energy between a
finite plate, of dimensions stated above, and an infinite plate is: 
\[
{\cal{E}}=2\sum_{m=1}^{m_0}{\cal{E}}_{2m}^1, 
\]
when $m_0\le m_1$ and
\[
{\cal{E}}=2\sum_{m=1}^{m_1-1}{\cal{E}}_{2m}^1+2{\cal{E}}_{2m_1}^0 
\]
when $m_0> m_1$, where
\[
{\cal{E}}_{2m}^1=-\frac{\hbar cW\cos^2(m\gamma)\sin\gamma}{64\pi^2\sin^4(m\gamma)}\left(\frac{1}{R_0^2\cos\gamma}-\frac{1}{R_1^2\cos[(m-1)\gamma]\cos(m\gamma)}\right),
\]
and
\[
{\cal{E}}_{2m_1}^0=\frac{\hbar
  cW\cos^2(m_1\gamma)(R_0\cos\gamma-R_1\cos[(m_1-1)\gamma])^2}{64\pi^2\sin^4(m_1\gamma)\cos\gamma\cos[(m_1-1)\gamma]\sin(m_1\gamma)R_0^2R_1^2},  
\]
for $m_1$ even, and
\[
{\cal{E}}_{2m_1}^0=\frac{\hbar
  cW\cos^2(m_1\gamma)(R_0-R_1\cos[(m_1-1)\gamma])^2}{64\pi^2\sin^4(m_1\gamma)\sin[(m_1-1)\gamma]\cos[(m_1-1)\gamma]R_0^2R_1^2},
\]
for $m_1$ odd.
}

\vspace{0.2in}

In particular, the Casimir force is always attractive. We also prove that: 

\vspace{0.2in}

\noindent {\bf Main Theorem 3}.

{\it
In the parallel plate limit the above 
energy reduces to Casimir's original expression:
\[
{\cal{E}}=-\frac{\hbar c\pi^2A}{1440\;L^3},
\]
where $A$ is the area of the finite plate and $L$ is the separation
between the plates. 
}

\vspace{0.2in}

We now sketch the steps of our calculation. In Propositions
\ref{p:evenlength} and \ref{p:oddlength}, we reduce the infinite sum
over closed paths in the Casimir energy to a finite sum and compute
the total length of a closed path as a function the opening angle and
initial point. This establishes Main Theorem 1.

By virtue of the general
reflection formulae derived in \cite{gak4} and summarised in
Proposition \ref{p:genref}, we compute the sequence of rays reflected  
$k$-times off the pair of planes. The results (Proposition
\ref{p:genit}) are then used to determine the initial direction
required for the ray to return to a given point in the wedge after 
$k$-bounces (Propositions \ref{p:evencl} and \ref{p:oddcl}). In
Propositions \ref{p:evenseq} and \ref{p:oddseq} we find the sequence
of intersection points in the wedge. With the aid of these, we can
determine the restrictions on the regions of integration in the energy
for a finite plate above an infinite plane. This
allows us to integrate explicitly and find the Casimir energy as a
function of the dimensions and relative positions of the plates, as
stated in Main Theorem 2.

This paper is organised as follows. In the next section we review the
optical approximation to the Casimir effect, as developed by Jaffe and
Scardicchio, and use the classical method of images to determine the
existence and lengths of closed paths in an infinite wedge.  

Section 3 outlines our approach to the geometry of
reflection in terms of the complex geometry on
$\mbox{T}{\Bbb{P}}^1$ - further details can be found in \cite{gak4}. 
It also outlines how our technique can be used to compute the
Casimir energy when the boundaries are not planar.
In Section 4 we use this geometry to find the sequence
of intersection points of a closed $n$-bounce path in a wedge. Here we
must treat even and odd bounce closed paths separately.

Sections 5 combines these results to determine the Dirichlet Casimir
energy between a plane and a finite plate. In Section 6 we obtain
Casimir's original result in the parallel plate limit (Main Theorem 3).

\vspace{0.2in} 

\noindent {\bf Acknowledgements}: The first two authors would like to
thank Karl Luttinger for many inspiring and helpful discussions on the
mathematical background of this paper. The proof of Main Theorem 1
was significantly shortened by suggestions of Patrick Dorey and
Michael Farber. 

The first author would like to express his
appreciation for the hospitality of Grey College, Durham University,
during the development of this work. This work was made possible by
the International Collaboration Programme, Enterprise Ireland.

\section{The Geometric Optics Approximation to the Casimir Effect in a
Wedge}

The Casimir effect is a force that appears between conducting
boundaries due to quantum fluctuations of the vacuum energy of the
electromagnetic field. Recently, Jaffe and Scardicchio \cite{jas} have
given the following explicit formula for the Casimir energy in terms
of geometric quantities:  
\[
{\cal{E}}=-\frac{\hbar c}{2\pi^2}\sum_n(-1)^n\iiint\limits_{{\cal{D}}_n}\frac{\sqrt{\Delta_n(x,x)}}{l_n^3(x)}d^3x.
\]
Here the sum is over straight paths with $n$-bounces between the boundaries which begin and
end at the point $x$, $l_n$ is the length of
such a path and $\Delta_n(x,x)$ is the expansion factor, or Van Vleck
determinant:
\[
\Delta_n(x,x')=\lim_{\delta\rightarrow 0}\delta^2e^{-2\int_\delta^lHdr},
\]
where $H$ is the mean curvature of the wavefront and $r$ is an affine
parameter along the path joining $x$ and $x'$. 

The computation of the Casimir effect can, in this approximation, be
reduced to determining the length $l_n$, Van Vleck determinant
$\Delta_n$ and region of existence ${\cal{D}}_n$ for closed
$n$-bounce paths between some prescribed boundary components.
When the boundaries are planes the Van Vleck determinant is just the
length of the path, so we need only compute $l_n$ and ${\cal{D}}_n$.

We now consider the Casimir effect in a wedge formed by two
non-parallel planes. Since any closed path in the wedge will be
contained in a plane perpendicular to the line of intersection of the
planes, the problem reduces to a 2-dimensional one. 

Consider, then, a wedge in the plane with opening angle $\gamma$, where
$0<\gamma<\frac{\pi}{2}$. The following two propositions determine the
existence and length of closed paths as a function of $\gamma$ and the
number of reflections, by applying the classical method of images.

\begin{Prop}\label{p:evenlength}
There exists a closed $2m$-bounce path in a wedge with angle
$\gamma$ iff $\gamma<\frac{\pi}{2m}$. There exists at most two closed paths 
(traversed in opposite directions) and the length of such a
closed path is $2R|\sin(m\gamma)|$. 
\end{Prop}
\begin{pf}
We reflect the wedge
$2m$-times through one of its sides. For a point ($R$, $\psi$) in the
wedge, we also reflect it $2m$-times. The image of its path reflected
in the wedge is a straight line joining the initial point to the final
image point. 

\vspace{0.2in}

\setlength{\epsfxsize}{3in}
\begin{center}
   \mbox{\epsfbox{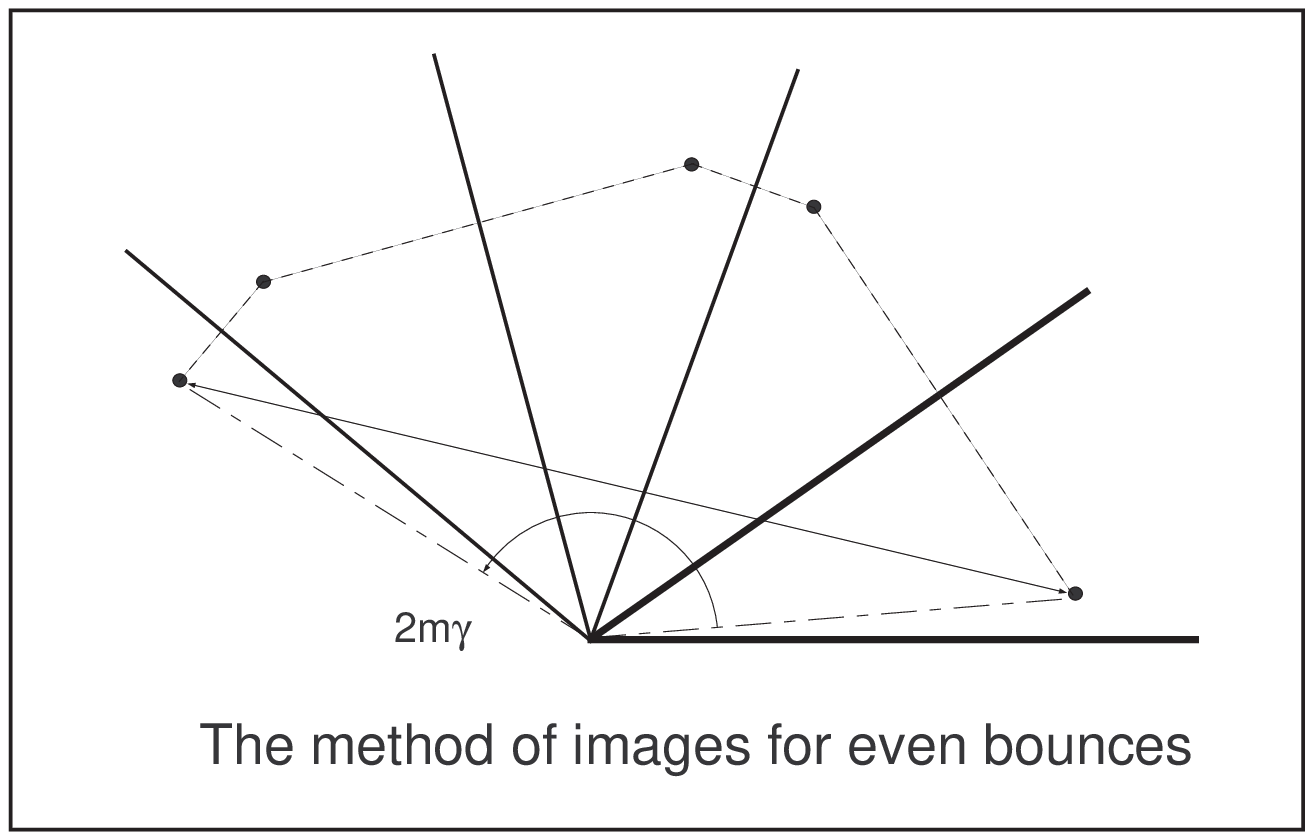}}
\end{center}

\vspace{0.2in}

The angle between the ray from the origin to the beginning and
endpoints is $2m\gamma$ and so the closed path exists iff
$2m\gamma<\pi$, which proves the first part of the
proposition. Applying the cosine rule we find that the distance
between these two points, which is equal to the length of the closed
path, is $2R|\sin(m\gamma)|$, as claimed. 
\end{pf}

The existence and length of a closed odd-bounce path is somewhat
different. Let ($R$, $\psi$) be polar coordinates in the plane, where
$R$ is the distance from the vertex, the bottom plane is aligned with
the horizontal and $\psi$ measures the angle from the vertical.

\begin{Prop}\label{p:oddlength}
There exists a closed $2m-1$-bounce path starting at the
point ($R\sin\psi  $, $R\cos\psi$) in a wedge of opening angle
$\gamma$ iff $\gamma<\frac{\pi}{2m-2}$ and either of the
following hold: $\psi<\pi-m\gamma$ or $\psi>(m-1)\gamma$.
The former applies to paths that first strike the non-horizontal plane,
while the latter applies to paths that first strike the horizontal plane.

There exists exactly none, one or two closed $2m-1$-bounce paths from a given point,
according to the inequalities above.

The total length of a closed $2m-1$ reflection path in a wedge with
angle $\gamma$ starting at the point ($R\sin\psi$, $R\cos\psi$) is
either $2R|\cos(\psi+m\gamma)|$  or $2R|\cos[\psi-(m-1)\gamma]|$
(respectively).  
\end{Prop}

\begin{pf}
Applying the classical method of images again, we reflect the wedge
$2m-1$-times through one of its sides. For a point ($R$, $\psi$) in
the wedge, we also reflect it $2m-1$-times in the non-horizontal plane
first. The image of its path reflected in the wedge is a straight line
joining the initial point to the final image point. 

\vspace{0.2in}

\setlength{\epsfxsize}{3in}
\begin{center}
   \mbox{\epsfbox{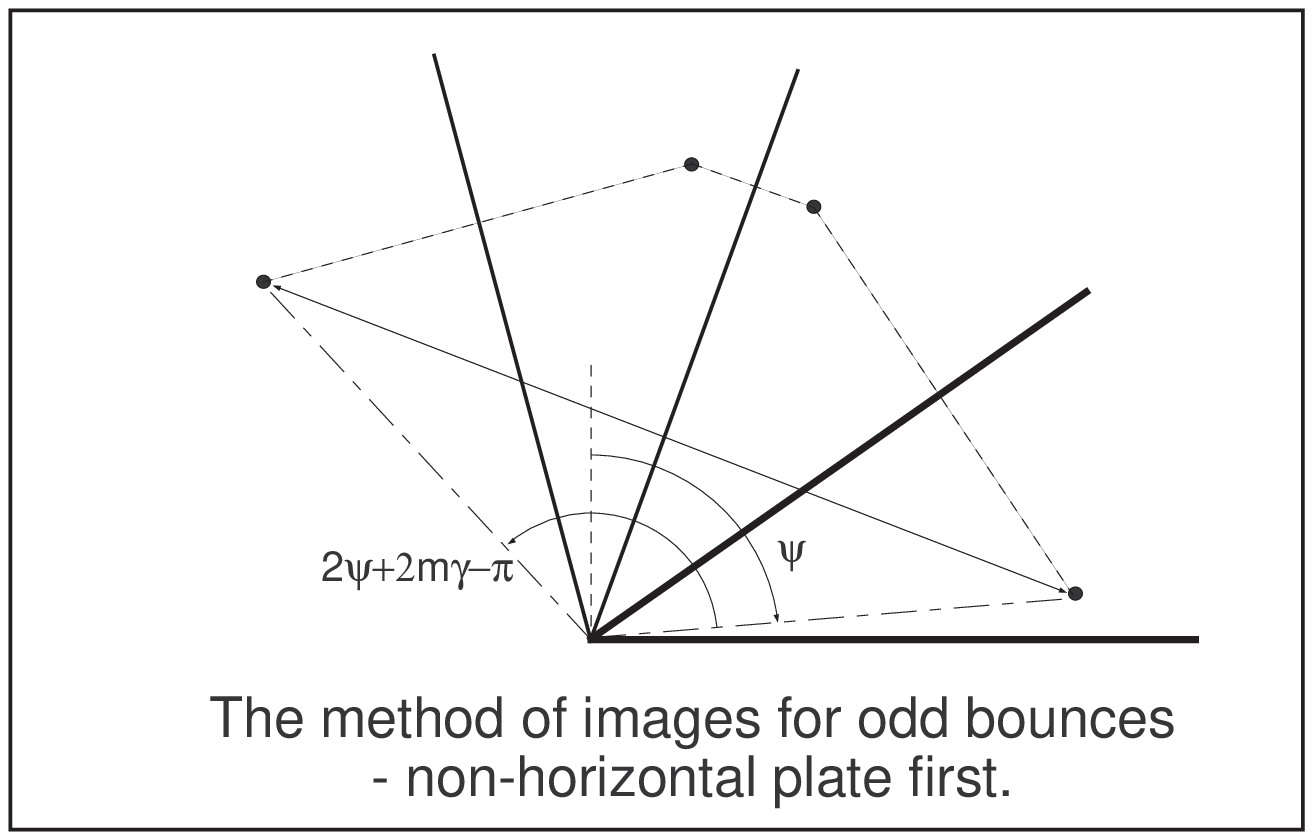}}
\end{center}

\vspace{0.2in}

The angle between the beginning and endpoints is
$2\psi-\pi+2m\gamma$ and so the closed path exists iff
$2\psi-\pi+2m\gamma<\pi$, which proves the first part of the
proposition. Applying the cosine rule we find that the
distance between the two points, which is equal to the length of the
closed path, is $2R|\cos(\psi+m\gamma)|$, as claimed. 

Similarly for the path that first strikes the horizontal plane.

\end{pf}

\begin{Cor}
A closed odd bounce path retraces itself.
\end{Cor}
\begin{pf}
It is clear from symmetry that a
reflected odd bounce path crosses the middle plane image at right
angles, and so is reflected back along itself.
\end{pf}

By way of example, the diagram below shows the regions where a closed
3-bounce exists for various opening angles $\gamma$.

\vspace{0.2in}

\setlength{\epsfxsize}{5in}
\begin{center}
   \mbox{\epsfbox{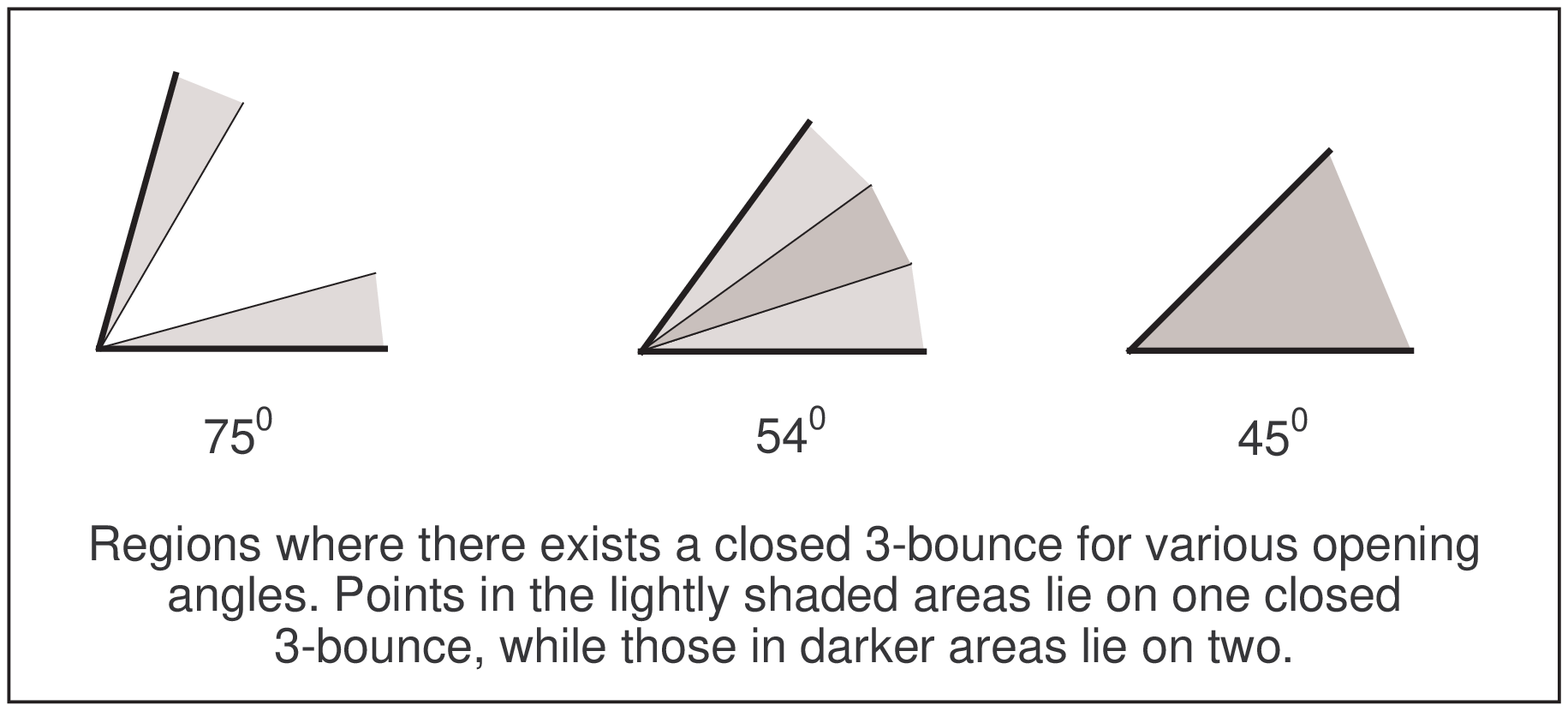}}
\end{center}

\vspace{0.2in}

Introducing the preceding geometric expression for the path lengths
into Jaffe and Scardicchio's optical approximation for the Casimir
energy, we find that for a wedge with opening angle $\gamma$ satisfying
$\frac{\pi}{2m+2}\le\gamma<\frac{\pi}{2m+2}$: 
\begin{align}
{\cal{E}}=&-\frac{\hbar c}{16\pi^2}\sum_{n=1}^m\quad\iiint\limits_{\frac{\pi}{2}-\gamma\le\psi\le \frac{\pi}{2}}\frac{1}{R^4\sin^4(n\gamma)}d^3x\nonumber\\
&\qquad\qquad+\frac{\hbar c}{32\pi^2}\sum_{n=1}^m\quad\iiint\limits_{\psi_1<\psi<\frac{\pi}{2}}\frac{1}{R^4\cos^4[\psi+(n+1)\gamma]}d^3x\nonumber\\
&\qquad\qquad\qquad+\frac{\hbar c}{32\pi^2}\sum_{n=1}^m\quad\iiint\limits_{\frac{\pi}{2}-\gamma<\psi<\psi_2}\frac{1}{R^4\cos^4(\psi-n\gamma)}d^3x,\nonumber
\end{align}
where $\psi_1=\mbox{max}[n\gamma,\frac{\pi}{2}-\gamma]$ and
$\psi_2=\mbox{min}[\pi-(n+1)\gamma,\frac{\pi}{2}]$. 

Here, we have used the fact that even bounces count twice as they can be
traversed in either direction, while odd bounces give two 
contributions - one from reflecting off the horizontal plane first and
one from reflecting off the top plane first. We denote
these three contributions by ${\cal{E}}_{2m}$, ${\cal{E}}_{2m+1}^H$ and
${\cal{E}}_{2m+1}^T$, respectively. 

The 1-bounce paths
have been excluded from the sum as the energy density diverges for such
paths, since their lengths go to zero as one approaches the boundary.

The difficulty with the above integrals is that they are divergent at
 the limit $R=0$.  This is caused by the intersection 
of the two planes where the closed path lengths go to zero. Indeed, we
cannot expect the geometric optics approximation to be accurate near
corners.

In order to remove this difficulty we put a finite separation between
the two plates. This will alter the regions of integration
${\cal{D}}_n$ and this we investigate after introducing some complex
geometry in the next section.

\section{The Geometry of Reflection}

We now describe geometric optics with the aid of the space of oriented
affine lines in ${\Bbb{R}}^3$.  Further details of this approach,
which we only summarise below, can be found in \cite{gak4} \cite{gak3}. 

The space of oriented affine lines in ${\Bbb{R}}^3$ can be identified
with the tangent bundle to the 2-sphere $\mbox{T}{\Bbb{P}}^1$
\cite{hitch}. From our point of view, geometric optics is the study of
2-parameter families of oriented lines or {\it line congruences}. These
lines, which we consider as a surface in the non-compact 4-manifold
$\mbox{T}{\Bbb{P}}^1$, form the rays of the optical system and the
wavefronts of the optical system are orthogonal to the line
congruence. 

Not every line congruence has such orthogonal wavefronts - the lines
may be twisting. Using the round metric on ${\Bbb{P}}^1$, the canonical
symplectic structure on $\mbox{T}^*{\Bbb{P}}^1$ can be pulled back to 
$\mbox{T}{\Bbb{P}}^1$. A line congruence is non-twisting iff this
symplectic structure vanishes on the associated surface
in $\mbox{T}{\Bbb{P}}^1$, i.e. it is Lagrangian.

In addition, $\mbox{T}{\Bbb{P}}^1$ carries a canonical complex
structure given by rotation through 90$^0$ about the line. This
preserves the tangent space to a line congruence at a given line 
iff the line is shearfree. In the case of a Lagrangian line
congruence, the associated wavefronts have an umbilical point along
this line.

We now consider reflection of a wavefront off a surface
in ${\Bbb{R}}^3$. First choose coordinates $\xi$ on ${\Bbb{P}}^1$ by
stereographic projection from the south pole of the unit sphere about
the origin onto the plane through the equator, and coordinates on
$\mbox{T}{\Bbb{P}}^1$ by identifying ($\xi$, $\eta$)$\in{\Bbb{C}}^2$ with 
\[
\eta\frac{\partial}{\partial \xi}+\overline{\eta}\frac{\partial}{\partial \overline{\xi}}\in \mbox{T}_\xi {\Bbb{P}}^1.
\]
Thus $\xi$ gives the direction of the oriented line in ${\Bbb{R}}^3$
and $\eta$ gives the perpendicular distance vector from the origin to
the line. Let ($x^1,x^2,x^3$) be Euclidean coordinates on
${\Bbb{R}}^3={\Bbb{C}}\oplus{\Bbb{R}}$, and  set $z=x^1+ix^2$,
$t=x^3$. The point ($z$, $t$) lies on the line ($\xi$, $\eta$) iff 
\begin{equation}\label{e:incidence}
\eta=\frac{1}{2}(z-2t\xi-\overline{z}\xi^2),
\end{equation}

Consider an incoming ray ($\xi_k$, $\eta_k$)$\in \mbox{T}{\Bbb{P}}^1$
reflecting off an oriented surface at a point ($\alpha_k$,
$b_k$)$\in{\Bbb{C}}\oplus{\Bbb{R}}={\Bbb{R}}^3$. Suppose that the
oriented normal to the surface at the point of reflection is ($\nu_k$,
$\chi_k$)$\in \mbox{T}{\Bbb{P}}^1$ and that the reflected ray is
($\xi_{k+1}$, $\eta_{k+1}$)$\in \mbox{T}{\Bbb{P}}^1$. We denote the
(oriented) distance of ($\alpha_k$, $b_k$) from the closest point to the
origin on the incoming ray, reflected ray and normal by $r_k$,
$r_{k+1}$ and $s_k$, respectively (see the diagram below).  

\vspace{0.2in}

\setlength{\epsfxsize}{3in}
\begin{center}
   \mbox{\epsfbox{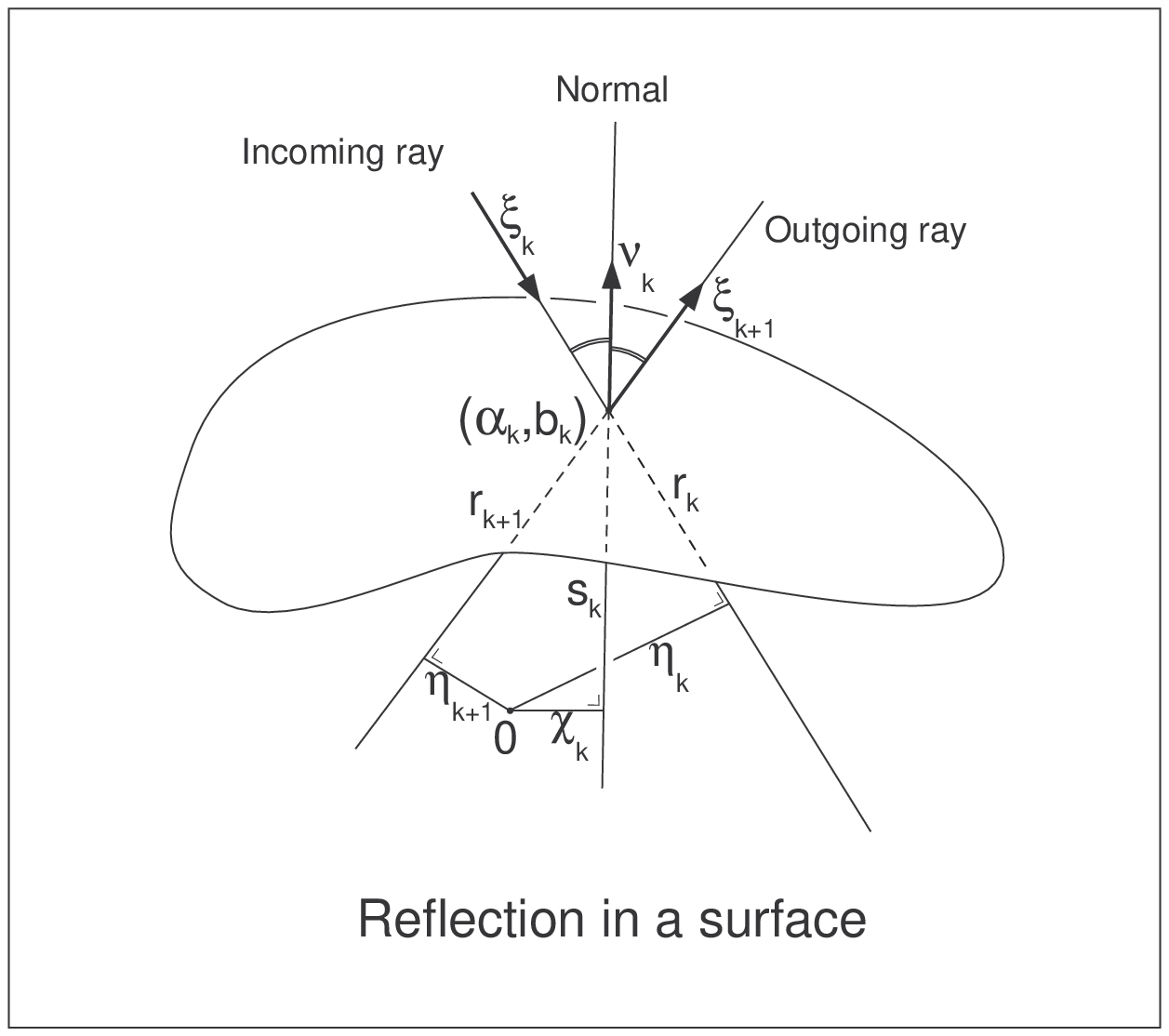}}
\end{center}

\vspace{0.2in}

The following proposition describes reflection as the combined actions of 
$\overline{\mbox{PSL}(2,{\Bbb{C}})}$ and fibre mappings on
$\mbox{T}{\Bbb{P}}^1$.

\begin{Prop}\label{p:genref}\cite{gak4}
The reflected ray is given by
\begin{equation}\label{e:refdir}
\xi_{k+1}=\frac{2\nu_k\bar{\xi}_k+1-\nu_k\bar{\nu}_k}
           {(1-\nu_k\bar{\nu}_k)\bar{\xi}_k-2\bar{\nu}_k},
\end{equation}
\begin{equation}\label{e:refperp}
\eta_{k+1}=\frac{-(1+\nu_k\bar{\nu}_k)^2\bar{\eta}_k
+(\bar{\nu}_k-\bar{\xi}_k)(1+\nu_k\bar{\xi}_k)(1+\nu_k\bar{\nu}_k)s_k
}{((1-\nu_k\bar{\nu}_k)\bar{\xi}_k-2\bar{\nu}_k)^2},
\end{equation}
the distance of ($\alpha_k$, $b_k$) from the closest point to the
origin on the reflected ray is
\begin{equation}\label{e:refpot}
r_{k+1}=r_k+\frac{2(|\nu_k-\xi_k|^2-|1+\nu_k\bar{\xi}_k|^2)}
    {(1+\nu_k\bar{\nu}_k)(1+\xi_k\bar{\xi}_k)}\;s_k,
\end{equation}
and the intersection equation is 
\begin{equation}\label{e:int}
\eta_k=\frac{(1+\bar{\nu}_k\xi_k)^2\chi_k
-(\nu_k-\xi_k)^2\bar{\chi}_k
+(\nu_k-\xi_k)(1+\bar{\nu}_k\xi_k)(1+\nu_k\bar{\nu}_k)s_k
}{(1+\nu_k\bar{\nu}_k)^2}.
\end{equation}
\end{Prop}

There are many equivalent ways of rewriting and using these
equations. In \cite{gak4} these equations are used to explicitly
determine the scattering of plane and spherical waves off planes,
spheres and tori. To do this we first solve the intersection
equation (\ref{e:int}) to find the normal $\nu_k$ to the surface at
the point of reflection and then use (\ref{e:refdir}) and
(\ref{e:refperp}) to determine the outgoing ray.

The technique above can be utilised to compute the Casimir
effect between  non-planar boundaries. To get an explicit general form for
the Van Vleck determinant in our formalism we proceed as follows. As
detailed in \cite{gak2}, the mean curvature of a wavefront orthogonal
to the line congruence 
$(\zeta,\bar{\zeta})\rightarrow(\xi(\zeta,\bar{\zeta}),\eta(\zeta,\bar{\zeta}))$ is
\[
H=\frac{ \partial^+\eta\;\bar{\partial}\;\bar{\xi} -\partial^-\eta\;\partial\bar{\xi}}
{\partial^+\eta\;\overline{\partial^+\eta}-\partial^-\eta\;\overline{\partial^-\eta}},
\]
where $\partial=\frac{\partial}{\partial\zeta}$ and 
\[
\partial^+\eta\equiv\partial \eta-\frac{2\eta\overline{\xi}\partial \xi}{1+\xi\overline{\xi}}+r\partial\xi,
\qquad\qquad
\partial^-\eta\equiv\overline{\partial} \eta-\frac{2\eta\overline{\xi}\;\overline{\partial} \xi}{1+\xi\overline{\xi}}+r\overline{\partial}\xi.
\]
The Lagrangian condition for the line congruence is equivalent to $H$
being real. 

Suppose we start with a spherical wave emanating from the point
($\alpha_0$, $b_0$)$\in{\Bbb{R}}^3$ and this is reflected $n$ times
before passing through a point ($\alpha_{n+1}$, $b_{n+1}$). Let
($\xi_k$, $\eta_k$) be the line congruence after the $k$th reflection.

\begin{Prop}
The geometric form of the Van Vleck determinant between the start and
end points is
\[
\Delta_n(0,n+1)=\frac{1}{l_0^2}\prod_{k=1}^n\frac{[\Psi_k]_1}{[\Psi_k]_2},
\]
where $l_0$ is the distance from the initial point to the first reflection,
\[
\Psi_k=\partial^+\eta\;\overline{\partial^+\eta}-\partial^-\eta\;\overline{\partial^-\eta},
\]
$[\Psi_k]_1$ is evaluation of $\Psi_k$ at the start of the $k^{th}$
reflection and $[\Psi_k]_2$ is evaluation of $\Psi_k$ at the end of
the $k^{th}$ reflection. 
\end{Prop}
\begin{pf}
This follows immediately from the fact that the mean curvature is the
derivative with respect to $r$ of the natural log of $\Psi_k$. The
term $\delta^{-2}$ is cancelled by the singularity in the spherical
wavefront at the initial point. 
\end{pf}

\begin{Cor}
The Van Vleck determinant of two points
between plane boundaries is the reciprocal of the path length squared.
\end{Cor}

\begin{pf}
A spherical wavefront reflected in a plane remains spherical with the
same principal curvatures at the point of reflection and so
$[\Psi_{k}]_2=[\Psi_{k+1}]_1$. Moreover,
$[\sqrt{\Psi_{k}}]_2=[\sqrt{\Psi_{k}}]_1+l_k$, where $l_k$ is the
length of the k$^{th}$ reflected path. Thus
\[
\Delta_n(0,n+1)=\frac{1}{[\Psi_n]_2}=\left(\sum_{k=0}^{k=n}l_k\right)^{-2
},
\]
as claimed.
\end{pf}

For the background of this aspect of geometric optics see
Chapter 5 and Appendix B of \cite{kak}.

\section{Reflections in an Infinite Wedge}

Consider a ray originating from ($\alpha_0$, $b_0$) and
reflecting off a series of planes. Suppose that the points of
reflection on the planes are given by ($\alpha_k$, $b_k$), the
reflected directions are $\xi_{k+1}$ and $l_{k\;k+1}$ are the
(oriented) distances between these points, for $k=1,2,...$. Let
$\nu_k$ be the normal direction of the $k$th plane.

\begin{Prop}
The points of reflection are given by
\begin{equation} \label{e:genab}
\alpha_k=\alpha_{k-1}+\frac{2\xi_k}{1+\xi_k\bar{\xi}_k}l_{k-1\;k},
\qquad\qquad
b_k=b_{k-1}+\frac{1-\xi_k\bar{\xi}_k}{1+\xi_k\bar{\xi}_k}l_{k-1\;k},
\end{equation}
where
\begin{equation} \label{e:genlength}
l_{k-1\;k}=-\frac{(1+\xi_k\bar{\xi}_k)(\alpha_{k-1}\bar{\nu}_k
 +\bar{\alpha}_{k-1}\nu_k+(1-\nu_k\bar{\nu}_k)b_{k-1}-(1+\nu_k\bar{\nu}_k)s_k)}
  {|1+\nu_k\bar{\xi}_k|^2-|\nu_k-\xi_k|^2}.
\end{equation}
\end{Prop}
\begin{pf}
The first two equations hold by the definition of $l_{k-1\;k}$ as being
the distance between, and $\xi_k$ the direction of the oriented line
joining,  $(\alpha_k,b_k)$ and $(\alpha_{k-1},b_{k-1})$.
The normal line contains the point $(\alpha_k , b_k)$ and the incoming ray 
contains the point  $(\alpha_{k-1} , b_{k-1})$. The corresponding
incidence relations (\ref{e:incidence}) are 
\[
\chi_k = \frac{1}{2} (\alpha_k -2 b_k\nu_k - \bar\alpha_k \nu_k^2), 
\qquad\qquad
\eta_k = \frac{1}{2} (\alpha_{k-1} -2 b_{k-1}\xi_k - \bar\alpha_{k-1} \xi_k^2 ).
\]

Introducing these relations together with the first two of the
Proposition into equation (\ref{e:int}) gives the third.

\end{pf}

From here on we consider only the case of two plane boundaries.
Consider now reflection in two planes, one of which is horizontal, 
both containing the origin.  Then $\nu_1=0$ and we let $\nu_2$ be the normal direction to the
non-horizontal plane pointing inward in the first quadrant. 

\begin{Prop}
The lengths of the reflected paths satisfy
\begin{equation}\label{e:oddlength}
l_{2k-1\;2k}=-\frac{(1+\xi_{2k}\bar{\xi}_{2k})(\alpha_{2k-1}\bar{\nu}_2
 +\bar{\alpha}_{2k-1}\nu_2)}
  {(1-\xi_{2k}\bar{\xi}_{2k})(1-\nu_{2}\bar{\nu}_{2})+2\xi_{2k}\bar{\nu}_2
    +2\bar{\xi}_{2k}\nu_2},
\end{equation}
\begin{equation}\label{e:evenlength}
l_{2k\;2k+1}=\frac{(1+\xi_{2k+1}\bar{\xi}_{2k+1})(1-\xi_{2k}\bar{\xi}_{2k})}   {(1-\xi_{2k+1}\bar{\xi}_{2k+1})(1+\xi_{2k}\bar{\xi}_{2k})}l_{2k-1\;2k}.
\end{equation}
\end{Prop}
\begin{pf}
 The first equation follows from equation (\ref{e:genlength}) and
the fact that $s_k = 0 $ and $b_{2k-1}= 0$. The second equation follows 
from the same equation using $\nu_{2k+1} = 0$ and equation
the second of (\ref{e:genab}) in equation (\ref{e:genlength}).
\end{pf}

In the case of two planes there is a translational symmetry which we now exploit.
We assume that the line of intersection of the
planes lies along the $x^1$-axis and the acute angle lies in the first
quadrant. Set $\nu_2=tan(\beta/2)$ so that $\gamma=\pi-\beta$ is the
opening angle of the wedge.

With these simplifications the problem reduces to reflection in two
lines in the $x^2x^3$-plane, and $\alpha_k$, $\xi_k$ and $\nu_k$ are
all real. Let $\alpha_k=a_k\in{\Bbb{R}}$ and introduce polar coordinates ($R$,
$\psi$) in the $x^2x^3$-plane.

The reflected rays after $k$-bounces are given by:

\begin{Prop}\label{p:genit}
Consider a ray with direction $\xi_1$ emanating from the point
($ a_0$, $b_0$) and striking the horizontal plane. The sequence of
reflected rays is given by
\begin{equation}\label{e:genray}
\xi_{2k}=\frac{-\sin[(k-1)\beta\xi_1 +\cos[(k-1)\beta]}
   {\cos[(k-1)\beta]\xi_1+\sin[(k-1)\beta] },
\qquad
\xi_{2k+1}=\frac{\cos[k\beta]\xi_1+\sin[k\beta] }
   {\sin[k\beta]\xi_1 +\cos[k\beta]},
\end{equation}
\begin{equation}\label{e:genperp}
\eta_{2k}=\frac{-(a_0-2b_0\xi_1- a_0\xi_1^2)}
  {2(\cos[(k-1)\beta]\xi_1+\sin[(k-1)\beta] )^2},
\quad
\eta_{2k+1}=\frac{ a_0-2b_0\xi_1-a_0\xi_1^2}
  {2(\sin[k\beta]\xi_1 +\cos[k\beta])^2}.
\end{equation}
\end{Prop}

\begin{pf}
The first of equations (\ref{e:genray}) and (\ref{e:genperp}) are true for $k=0$.
We now proceed inductively and assume they both hold for
$k=k_0$. Then, since  
$\nu_{2k+1}=0$ and $\nu_{2k}=\nu_2$, we have by (\ref{e:refdir}) that
\[
\xi_{2k_0+2}=\frac{1}{\xi_{2k_0+1}}=\frac{-\sin(k_0\beta)\bar{\xi}_1+\cos(k_0\beta)}
   {\cos(k_0\beta)\bar{\xi}_1+\sin(k_0\beta)},
\]
which is the first of (\ref{e:genray}) for $k=k_0$. Similarily,
\begin{align}
\xi_{2k_0+3}=&\frac{2\nu_2\bar{\xi}_{2k_0+2}+1-\nu_2\bar{\nu}_2}
           {(1-\nu_2\bar{\nu}_2)\bar{\xi}_{2k_0+2}-2\bar{\nu}_2}\nonumber\\
&=\frac{\sin\beta[-\sin(k_0\beta)\xi_1+\cos(k_0\beta)]+\cos\beta[\cos(k_0\beta)\xi_1+\sin(k_0\beta)]}
   {\cos\beta[-\sin(k_0\beta)\xi_1+\cos(k_0\beta)]-\sin\beta[\cos(k_0\beta)\xi_1+\sin(k_0\beta)]}\nonumber\\&=\frac{[\cos\beta\cos(k_0\beta)-\sin\beta\sin(k_0\beta)]\xi_1+\sin\beta\cos(k_0\beta)+\cos\beta\sin(k_0\beta)}
{[-\cos\beta\sin(k_0\beta)-\sin\beta\cos(k_0\beta)]\xi_1+\cos\beta\cos(k_0\beta)-\sin\beta\sin(k_0\beta)}\nonumber\\
&=\frac{\cos[(k_0+1)\beta]\xi_1+\sin[(k_0+1)\beta]}
    {-\sin[(k_0+1)\beta]\xi_1+\cos[(k_0+1)\beta]}\nonumber,
\end{align}
which is the first of equation (\ref{e:genray}) with $k=k_0+1$. Thus
(\ref{e:genray}) hold for all $k$.

A similar inductive argument shows that, by (\ref{e:refperp}),
equations (\ref{e:genperp}) hold for all
$k$. 

\end{pf}

For future calculations we note that:

\begin{Lem}
If $\xi_1=tan(\phi/2) $ the sequence of reflected directions
satisfy 
\begin{equation}\label{e:evencossin}
\frac{1-\xi_{2k}\bar{\xi}_{2k}}{1+\xi_{2k}\bar{\xi}_{2k}}
     =-\cos[\phi+2(k-1)\beta], \qquad \frac{2\xi_{2k}}{1+\xi_{2k}\bar{\xi}_{2k}}
     =\sin[\phi+2(k-1)\beta],
\end{equation}
\begin{equation}\label{e:oddcossin}
\frac{1-\xi_{2k+1}\bar{\xi}_{2k+1}}{1+\xi_{2k+1}\bar{\xi}_{2k+1}}
     =\cos[\phi+2k\beta], \qquad\frac{2\xi_{2k+1}}{1+\xi_{2k+1}\bar{\xi}_{2k+1}}
     =\sin[\phi+2k\beta]. 
\end{equation}
\end{Lem}
\begin{pf}
These follow, with the aid of trigonometry identities, 
from equations (\ref{e:genray}).
\end{pf}
\vspace{0.2in}

\subsection{Even Reflections}

We now study paths that return to the original point after $2m$
reflections. The direction of the initial ray is given by:

\begin{Prop}\label{p:evencl}
A point ($ a_0$, $b_0$) lies on a closed path with $2m$
reflections iff the initial direction of the ray is
\begin{equation}\label{e:evenclosedray1}
\xi_1=\frac{\sin[\psi-m\beta]\pm 1}{\cos[\psi-m\beta]}, 
\end{equation}
where $ a_0=R\sin\psi  $ and $b_0=R\cos\psi$. For
$\xi_1=tan(\phi/2) $ this is equivalent to
\begin{equation}\label{e:evenclosedray2}
\cos\phi=\mp\sin[\psi-m\beta], \qquad\qquad \sin\phi=\pm\cos[\psi-m\beta].
\end{equation}
\end{Prop}
\begin{pf}

The initial point ($ a_0$, $b_0$) is on a closed $2m$-bounce if it is 
contained on the final outgoing ray:
\[
\eta_{2m+1} = \frac{1}{2} (a_0 -2 b_0\xi_{2m+1} -a_0 \xi_{2m+1}^2 ).
\]
Substituting the first of equations (\ref{e:genray}) and
(\ref{e:genperp}) in this gives the quadratic equation
\[
\cos[\psi-m\beta]\xi_1^2-2\sin[\psi-m\beta]\xi_1-\cos[\psi-m\beta]=0.
\]
The solution to this is (\ref{e:evenclosedray1}), or, equivalently, 
(\ref{e:evenclosedray2}).
\end{pf}

For future reference we note that:
\begin{Lem}
For a closed 2$m$ reflection path with first reflection off the
horizontal plane
\begin{equation}\label{e:closedevencossin}
\frac{1-\xi_{2k}\bar{\xi}_{2k}}{1+\xi_{2k}\bar{\xi}_{2k}}
     =\pm\sin[\psi-(m-2k+2)\beta],\qquad\frac{2\xi_{2k}}{1+\xi_{2k}\bar{\xi}_{2k}}
     =\pm\cos[\psi-(m-2k+2)\beta], 
\end{equation}
\begin{equation}\label{e:closedoddcossin}
\frac{1-\xi_{2k+1}\bar{\xi}_{2k+1}}{1+\xi_{2k+1}\bar{\xi}_{2k+1}}
     =\mp\sin[\psi-(m-2k)\beta],\qquad\frac{2\xi_{2k+1}}{1+\xi_{2k+1}\bar{\xi}_{2k+1}}
     =\pm\cos[\psi-(m-2k)\beta]. 
\end{equation}
\end{Lem}
\begin{pf}
These follow from substituting (\ref{e:evenclosedray2}) in (\ref{e:evencossin}) and 
(\ref{e:oddcossin}).
\end{pf}

For a closed $2m$-bounce path, the sequence of points of intersection with
the boundaries and the length of the paths are given by:

\begin{Prop}\label{p:evenseq}
For a closed $2m$ reflection path with first reflection off the
horizontal plane the sequence of points of reflection and path
lengths are
\begin{equation}\label{e:aevencl}
 a_{2k}=\frac{R\cos(m\beta)\cos\beta}{\sin[\psi-(m-2k+1)\beta]},
\qquad\qquad
 a_{2k-1}=\frac{R\cos(m\beta)}{\sin[\psi-(m-2k+2)\beta]},
\end{equation}
\begin{equation}\label{e:bevenevencl}
b_{2k}=  -\frac{R\cos(m\beta)\sin\beta}{\sin[\psi-(m-2k+1)\beta]},\qquad\qquad b_{2k-1}=0,
\end{equation}
\begin{equation}\label{e:lengthevencl}
l_{k\;k+1}=\mp\frac{R\cos(m\beta)\sin\beta}{\sin[\psi-(m-k)\beta]\sin[\psi-(m-k+1)\beta]},
\end{equation}
\begin{equation}\label{e:sfevenlength}
l_{01}=\pm\frac{R\cos\psi}{\sin[\psi-m\beta]},
\qquad
l_{2m\;0}=\mp\frac{R\cos[\psi-\beta]}{\sin[\psi+(m-1)\beta]},
\end{equation}
where the signs are chosen to make the lengths positive.
\end{Prop}
\begin{pf}
The second of (\ref{e:bevenevencl}) is true because odd reflections
are off the horizontal plane and so have zero $x^3$-coordinate.

The distance from the initial point $(a_0,b_0)$ to the first reflection is
\[
l_{01}=-\frac{R\cos\psi}{\cos\phi}=\pm\frac{R\cos\psi}{\sin[\psi-m\beta]},
\]
using equation (\ref{e:evenclosedray2}). This proves the first of
(\ref{e:sfevenlength}). Then, by the first of (\ref{e:genab}),
\[
a_1=a_0+\sin\phi \;l_{01}=\frac{R(\cos\phi\sin\psi-\sin\phi\cos\psi)}{\cos\phi},
\]
which reduces to the second of (\ref{e:aevencl}) with $k=1$. Next, by
(\ref{e:oddlength}), (\ref{e:evencossin}) and (\ref{e:evenclosedray2}),
\begin{align}
l_{12}&=-\frac{(1+\xi_2\bar{\xi}_{2})( \alpha_1\bar{\nu}_2 +\bar{\alpha}_1\nu_2)}
  {(1-\xi_2\bar{\xi}_2)(1-\nu_2\bar{\nu}_{2})+2\xi_2\bar{\nu}_2
    +2\bar{\xi}_2\nu_2}\nonumber\\
&=\mp\frac{R\cos(m\beta)\sin\beta}
{\sin(\psi-m\beta)[\sin(\psi-m\beta)\cos\beta+\cos(\psi-m\beta)\sin\beta]}\nonumber\\
&=\mp\frac{R\cos(m\beta)\sin\beta}{\sin(\psi-m\beta)\sin[\psi-(m-1)\beta]}\nonumber,
\end{align}
which is (\ref{e:lengthevencl}) with $k=1$. Continuing on, we have
from the first of (\ref{e:genab})
\begin{align}
a_2&=a_1+\frac{2\xi_2}{1+\xi_2\bar{\xi}_2} l_{12}\nonumber\\
&=\frac{R\cos(m\beta)}{\sin(\psi-m\beta)}-\frac{\cos(\psi-m\beta)R\cos(m\beta)\sin\beta}
   {\sin(\psi-m\beta)\sin[\psi-(m-1)\beta]}\nonumber\\
&=\frac{R\cos(m\beta)[\sin[\psi-(m-1)\beta]-\cos(\psi-m\beta)R\cos(m\beta)\sin\beta]}{\sin(\psi-m\beta)\sin[\psi-(m-1)\beta]}\nonumber\\
&=\frac{R\cos(m\beta)\cos\beta}{\sin[\psi-(m-1)\beta]}\nonumber,
\end{align}
which is the first of (\ref{e:aevencl}) with $k=1$. Finally, by
(\ref{e:evenlength}), together with (\ref{e:evencossin}),
(\ref{e:oddcossin}) and (\ref{e:evenclosedray2}) 
\begin{align}
l_{23}&=-\frac{1+\xi_3\bar{\xi}_3}{1-\xi_3\bar{\xi}_3}
    \frac{1-\xi_2\bar{\xi}_2}{1+\xi_2\bar{\xi}_2}l_{12}\nonumber\\
&=\mp\frac{\sin(\psi-m\beta)}{\sin[\psi-(m-2)\beta]}
    \frac{R\cos(m\beta)\sin\beta}{\sin(\psi-m\beta)\sin[\psi-(m-1)\beta]}\nonumber\\
&=\mp\frac{R\cos(m\beta)\sin\beta}{\sin[\psi-(m-2)\beta]\sin[\psi-(m-1)\beta]}\nonumber,
\end{align}
which is (\ref{e:lengthevencl}) with $k=2$.

We now proceed inductively, assuming (\ref{e:aevencl}) and (\ref{e:bevenevencl}) hold for $k=k_0$, 
and (\ref{e:lengthevencl}) holds for $k=2k_0$.
First, by the first of (\ref{e:genab}), (\ref{e:evenclosedray2}) and
(\ref{e:closedoddcossin}) 
\begin{align}
a_{2k_0+1}&=a_{2k_0}+\frac{2\xi_{2k_0+1}}{1+\xi_{2k_0+1}\bar{\xi}_{2k_0+1}}l_{2k_0\;2k_0+1}\nonumber\\
&=\frac{R\cos(m\beta)\cos\beta}{\sin[\psi-(m-2k_0+1)\beta]}
  -\frac{\cos[\psi-(m-2k_0)\beta]R\cos(m\beta)\sin\beta}
         {\sin[\psi-(m-2k_0)\beta]\sin[\psi-(m-2k_0+1)\beta]}\nonumber\\
&=\frac{R\cos(m\beta)[\cos\beta\sin[\psi-(m-2k_0)\beta]-\cos[\psi-(m-2k_0)\beta]\sin\beta]}
         {\sin[\psi-(m-2k_0)\beta]\sin[\psi-(m-2k_0+1)\beta]}\nonumber\\
&=\frac{R\cos(m\beta)}
         {\sin[\psi-(m-2k_0)\beta]}\nonumber,
\end{align}
which proves the second of (\ref{e:aevencl}) with $k=k_0+1$. Next we have from
(\ref{e:oddlength}), (\ref{e:evenclosedray2}) and
(\ref{e:closedevencossin}) that
\begin{align}
l_{2k_0+1\;2k_0+2}&=-\frac{(1+\xi_{2k_0+2}\bar{\xi}_{2k_0+2})( a_{2k_0+1}\bar{\nu}_2 +a_{2k_0+1}\nu_2)}
  {(1-\xi_{2k_0+2}\bar{\xi}_{2k_0+2})(1-\nu_2\bar{\nu}_{2})+2\xi_{2k_0+2}\bar{\nu}_2
    +2\bar{\xi}_{2k_0+2}\nu_2}\nonumber\\
&=\mp\frac{R\cos(m\beta)\sin\beta}
 {\sin[\psi-(m-2k_0)\beta][\sin[\psi-(m-2k_0)\beta]\cos\beta+\cos[\psi-(m-2k_0)\beta]\sin\beta]}\nonumber\\
&=\mp\frac{R\cos(m\beta)\sin\beta}
         {\sin[\psi-(m-2k_0)\beta]\sin[\psi-(m-2k_0-1)\beta]}\nonumber,
\end{align}
which is (\ref{e:lengthevencl}) with $k=2k_0+1$. Again, by
the first of (\ref{e:genab}), (\ref{e:evenclosedray2}) and
(\ref{e:closedevencossin}) 
\begin{align}
a_{2k_0+2}&=a_{2k_0+1}+\frac{2\xi_{2k_0+2}}{1+\xi_{2k_0+2}\bar{\xi}_{2k_0+2}}l_{2k_0+1\;2k_0+2}\nonumber\\
&=\frac{R\cos(m\beta)}{\sin[\psi-(m-2k_0)\beta]}
-\frac{\cos[\psi-(m-2k_0)\beta]R\cos(m\beta)\sin\beta}{\sin[\psi-(m-2k_0)\beta]\sin[\psi-(m-2k_0-1)\beta]}\nonumber\\
&=\frac{R\cos(m\beta)[\sin[\psi-(m-2k_0-1)\beta]-\cos[\psi-(m-2k_0)\beta]\sin\beta]}
        {\sin[\psi-(m-2k_0)\beta]\sin[\psi-(m-2k_0-1)\beta]}\nonumber\\
&=\frac{R\cos(m\beta)\cos\beta}{\sin[\psi-(m-2k_0-1)\beta]}\nonumber,
\end{align}
which is the first of (\ref{e:aevencl}) with $k=k_0+1$. Finally, the first of
(\ref{e:bevenevencl}) follows simply from the second of (\ref{e:genab}) and
(\ref{e:lengthevencl}). 

The last equation follows from (\ref{e:bevenevencl}) with $k=m$ and
\[
l_{2m\;0}=\mp\frac{1+\xi_{2m+1}\bar{\xi}_{2m+1}}{1-\xi_{2m+1}\bar{\xi}_{2m+1}}(b_0-b_{2m}).
\]

\end{pf}

\subsection{Odd Reflections}

The direction of an initial ray that returns to the original point after $2m-1$
reflections is given by:

\begin{Prop}\label{p:oddcl}
A point $(a_0=R\sin\psi, b_0=R\cos\psi)$
is on a closed $2m-1$-bounce which strikes the horizontal plane first if
\begin{equation}\label{e:oddclosedray1}
\xi_1=\frac{\cos[(m-1)\beta]\pm 1}{\sin[(m-1)\beta]} ,
\end{equation}
or, equivalently,
\begin{equation}\label{e:oddclosedray2}
\sin\phi=\pm\sin[(m-1)\beta], \qquad\qquad \cos\phi=\mp\cos[(m-1)\beta].
\end{equation}
Note that the two solutions are antipodal.
\end{Prop}
\begin{pf}

The initial point ($ a_0$, $b_0$) is on a closed $2m$-bounce if it is 
contained on the final outgoing ray ({\it cf.} equation
(\ref{e:incidence})): 
\[
\eta_{2m} = \frac{1}{2} (a_0 -2 b_0\xi_{2m} -a_0 \xi_{2m}^2 ).
\]
Substituting the first of equations (\ref{e:genray}) and
(\ref{e:genperp}) in this gives the quadratic equation
\begin{align}
&[\sin\psi(-1+\cos[2(m-1)\beta])+\cos\psi\sin[2(m-1)\beta]]\bar{\xi}_1^2\nonumber\\
&\qquad+2[\cos\psi(-1-\cos[2(m-1)\beta])+\sin\psi\sin[2(m-1)\beta]]\bar{\xi}_1\nonumber\\
&\qquad-\sin\psi(-1+\cos[2(m-1)\beta])-\cos\psi\sin[2(m-1)\beta]=0\nonumber.
\end{align}
The solution to this is (\ref{e:evenclosedray1}), or, equivalently, 
(\ref{e:evenclosedray2}).
\end{pf}

For future use we note the following:
\begin{Lem}
For a closed $2m-1$ reflection path with first reflection off the
horizontal plane
\begin{equation}\label{e:closedoddevencossin}
\frac{1-\xi_{2k}\bar{\xi}_{2k}}{1+\xi_{2k}\bar{\xi}_{2k}}
     =\pm\cos[(m-2k+1)\beta],\qquad\frac{2\xi_{2k}}{1+\xi_{2k}\bar{\xi}_{2k}}
     =\pm\sin[(m-2k+1)\beta],
\end{equation}
\begin{equation}\label{e:closedoddoddcossin}
\frac{1-\xi_{2k+1}\bar{\xi}_{2k+1}}{1+\xi_{2k+1}\bar{\xi}_{2k+1}}
     =\mp\cos[(m-2k-1)\beta],\qquad\frac{2\xi_{2k+1}}{1+\xi_{2k+1}\bar{\xi}_{2k+1}}
     =\pm\cos[(m-2k-1)\beta] .
\end{equation}
\end{Lem}
\begin{pf}
These follow from substituting (\ref{e:oddclosedray2}) in (\ref{e:evencossin}) and 
(\ref{e:oddcossin}).
\end{pf}

For a closed $2m-1$-bounce path, the sequence of points of intersection with
the boundaries and the length of the paths are given by:

\begin{Prop}\label{p:oddseq}
For a closed $2m-1$-bounce with first reflection off the
horizontal plane the sequence of points of reflection and path
lengths are
\begin{equation}\label{e:aoddcl}
 a_{2k}=\frac{R\sin[\psi+(m-1)\beta]\cos\beta}{\cos[(m-2k)\beta]},
\qquad\qquad
 a_{2k-1}=\frac{R\sin[\psi+(m-1)\beta]}{\cos[(m-2k+1)\beta]},
\end{equation}
\begin{equation}\label{e:bevenoddcl}
b_{2k}= -\frac{R\sin[\psi+(m-1)\beta]\sin\beta}{\cos[(m-2k)\beta]}, \qquad\qquad b_{2k-1}=0,
\end{equation}
\begin{equation}\label{e:lengthoddcl}
l_{k\;k+1}=\mp\frac{R\sin[\psi+(m-1)\beta]\sin\beta}{\cos[(m-k)\beta]\cos[(m-k-1)\beta]},
\qquad\qquad
l_{01}=\pm\frac{R\cos\psi}{\sin[\psi-m\beta]},
\end{equation}
where the signs are chosen to make the lengths positive.
\end{Prop}

\begin{pf}
The second equation of (\ref{e:bevenoddcl}) holds since odd
reflections are off the horizontal plate, so the $x^3$-coordinate of
the intersection point is zero.

The distance from the initial point $(a_0,b_0)$ to the first reflection is
\[
l_{01}=-\frac{R\cos\psi}{\cos\phi}=\pm\frac{R\cos\psi}{\cos[(m-1)\beta]},
\]
using equation (\ref{e:oddclosedray2}). This proves
the second of (\ref{e:lengthoddcl}). Then, by the first of
(\ref{e:genab}) and (\ref{e:oddclosedray2}) 
\[
a_1=a_0+\sin\phi \;l_{01}=\frac{R\sin[\psi+(m-1)\beta]}{\cos[(m-1)\beta]},
\]
which is the second of (\ref{e:aoddcl}) with $k=1$. Next, by (\ref{e:oddlength}),
(\ref{e:oddclosedray2}) and (\ref{e:closedoddevencossin})
\begin{align}
l_{12}&=-\frac{(1+\xi_2\bar{\xi}_{2})( a_1\bar{\nu}_2 +a_1\nu_2)}
  {(1-\xi_2\bar{\xi}_2)(1-\nu_2\bar{\nu}_{2})+2\xi_2\bar{\nu}_2
    +2\bar{\xi}_2\nu_2}\nonumber\\
&=\mp\frac{R\sin[\psi+(m-1)\beta]\sin\beta}
{\cos[(m-1)\beta][\cos[(m-1)\beta]\cos\beta+\sin[(m-1)\beta]\sin\beta]}\nonumber\\
&=\mp\frac{R\sin[\psi+(m-1)\beta]\sin\beta}{\cos[(m-1)\beta]\cos[(m-2)\beta]}\nonumber,
\end{align}
which is the first of (\ref{e:lengthoddcl}) with $k=1$. Continuing on,
we have by the first of (\ref{e:genab}), (\ref{e:oddclosedray2}) and
(\ref{e:closedoddevencossin}) 
\begin{align}
a_2&=a_1+\frac{2\xi_2}{1+\xi_2\bar{\xi}_2} l_{12}\nonumber\\
&=\frac{R\sin[\psi+(m-1)\beta][\cos[(m-2)\beta]-\sin[(m-1)\beta]\sin\beta]}{\cos[(m-1)\beta]\cos[(m-2)\beta]}\nonumber\\
&=\frac{R\sin[\psi+(m-1)\beta]\cos[(m-1)\beta]\cos\beta}{\cos[(m-1)\beta]\cos[(m-2)\beta]}\nonumber\\
&=\frac{R\sin[\psi+(m-1)\beta]\cos\beta}{\cos[(m-2)\beta]}\nonumber,
\end{align}
which is the first of (\ref{e:aoddcl}) with $k=1$. Finally, by (\ref{e:evenlength}),
(\ref{e:closedoddevencossin}) and (\ref{e:closedoddoddcossin})
\begin{align}
l_{23}&=-\frac{1+\xi_3\bar{\xi}_3}{1-\xi_3\bar{\xi}_3}
    \frac{1-\xi_2\bar{\xi}_2}{1+\xi_2\bar{\xi}_2}l_{12}\nonumber\\
&=\mp\frac{\cos[(m-1)\beta]}{\cos[(m-3)\beta]}
    \frac{R\sin[\psi+(m-1)\beta]\sin\beta}{\cos[(m-1)\beta]\cos[(m-2)\beta]}\nonumber\\
&=\mp\frac{R\sin[\psi+(m-1)\beta]\sin\beta}{\cos[(m-3)\beta]\cos[(m-2)\beta]}\nonumber,
\end{align}
which is the first of (\ref{e:lengthoddcl}) with $k=2$.

We proceed inductively, assuming (\ref{e:aoddcl})
and (\ref{e:bevenoddcl}) hold for $k=k_0$, 
and (\ref{e:lengthoddcl}) holds for $k=2k_0$.
First, by the first of (\ref{e:genab}), (\ref{e:oddclosedray2}) and
(\ref{e:closedoddoddcossin}) 
\begin{align}
a_{2k_0+1}&=a_{2k_0}+\frac{2\xi_{2k_0+1}}{1+\xi_{2k_0+1}\bar{\xi}_{2k_0+1}}l_{2k_0\;2k_0+1}\nonumber\\
&=\frac{R\sin[\psi+(m-1)\beta]\cos\beta}{\cos[(m-2k_0)\beta]}
  -\frac{\sin[(m-2k_0-1)\beta]R\sin[\psi+(m-1)\beta]\sin\beta}
         {\cos[(m-2k_0)\beta]\cos[(m-2k_0-1)\beta]}\nonumber\\
&=\frac{R\sin[\psi+(m-1)\beta][\cos\beta\cos[(m-2k_0-1)\beta]-\sin[(m-2k_0-1)\beta]\sin\beta]}{\cos[(m-2k_0)\beta]\cos[(m-2k_0-1)\beta]}\nonumber\\
&=\frac{R\sin[\psi+(-1)m\beta]}{\cos[(m-2k_0-1)\beta]}\nonumber,
\end{align}
which proves the second of (\ref{e:aoddcl}) holds with $k=k_0+1$. Next, by
(\ref{e:oddlength}), (\ref{e:oddclosedray2}) and
(\ref{e:closedoddevencossin}), we have
\begin{align}
l_{2k_0+1\;2k_0+2}&=-\frac{(1+\xi_{2k_0+2}\bar{\xi}_{2k_0+2})( a_{2k_0+1}\bar{\nu}_2 +a_{2k_0+1}\nu_2)}
  {(1-\xi_{2k_0+2}\bar{\xi}_{2k_0+2})(1-\nu_2\bar{\nu}_{2})+2\xi_{2k_0+2}\bar{\nu}_2
    +2\bar{\xi}_{2k_0+2}\nu_2}\nonumber\\
&=\mp\frac{R\sin[\psi+(m-1)\beta]\sin\beta}
  {\cos[(m-2k_0-1)\beta][\cos[(2k_0-m+1)\beta]\cos\beta-\sin[(2k_0-m+1)\beta]\sin\beta]}\nonumber\\
&=\mp\frac{R\sin[\psi+(m-1)\beta]\sin\beta}
  {\cos[(m-2k_0-1)\beta]\cos[(m-2k_0-2)\beta]}\nonumber,
\end{align}
which is the first of (\ref{e:lengthoddcl}) with $k=2k_0+1$. Now by
the first of (\ref{e:genab}) and (\ref{e:closedoddevencossin})
\begin{align}
a_{2k_0+2}&=a_{2k_0+1}+\frac{2\xi_{2k_0+2}}{1+\xi_{2k_0+2}\bar{\xi}_{2k_0+2}}l_{2k_0+1\;2k_0+2}\nonumber\\
&=\frac{R\sin[\psi+(m-1)\beta]}{\cos[(m-2k_0-1)\beta]}
  -\frac{\sin[(m-2k_0-1)\beta]R\sin[\psi+(m-1)\beta]\sin\beta}
            {\cos[(m-2k_0-1)\beta]\cos[(m-2k_0-2)\beta]}\nonumber\\
&=\frac{R\sin[\psi+(m-1)\beta][\cos[(m-2k_0-2)\beta]
           -\sin[(m-2k_0-1)\beta]\sin\beta]}{\cos[(m-2k_0-1)\beta]\cos[(m-2k_0-2)\beta]}\nonumber\\
&=\frac{R\sin[\psi+(m-1)\beta]\cos\beta}{\cos[(m-2k_0-2)\beta]}\nonumber,
\end{align}
which is the first of (\ref{e:aoddcl}) with $k=k_0+1$. Finally, (\ref{e:bevenoddcl}) 
follows simply from the second of (\ref{e:genab}) and (\ref{e:lengthoddcl}).
\end{pf}

As proved earlier, a closed odd bounce retraces itself.
This can also be seen from the fact that for a closed $2m-1$-bounce, with $m$ odd, 
$ a_{2k-1}= a_{2(m-k)+1}$, while for $m$ even $ a_{2k}= a_{2(m-k)}$.

\begin{Prop}\label{p:oddseqt}
For a closed $2m-1$-bounce with first reflection off the
horizontal plane the sequence of points of reflection and path
lengths are
\begin{equation}\label{e:aoddclt}
 a_{2k}=\frac{R\sin[\psi-m\beta]}{\cos[(m-2k)\beta]},
\qquad\qquad
 a_{2k-1}=\frac{R\sin[\psi-m\beta]\cos\beta}{\cos[(m-2k+1)\beta]},
\end{equation}
\begin{equation}\label{e:bevenoddclt}
b_{2k-1}= -\frac{R\sin[\psi-m\beta]\sin\beta}{\cos[(m-2k+1)\beta]}, \qquad\qquad b_{2k}=0,
\end{equation}
\end{Prop}
\begin{pf}
These follow from (\ref{e:aoddcl}) and
(\ref{e:bevenoddcl}) by reflecting through the bisector of the wedge:
$(R,\psi)\rightarrow (R,\beta-\psi)$, so that 
\[
(a_{2k},b_{2k})\rightarrow\left(-\frac{a_{2k}}{\cos\beta},0\right),
\qquad\qquad
(a_{2k-1},0)\rightarrow\left(-a_{2k-1}\cos\beta,a_{2k-1}\sin\beta\right).
\]
\end{pf}

\vspace{0.2in}

\noindent{\bf Note}:

Instead of using the method of images as in Section 2, the lengths of the
closed $2m$- and $2m-1$-bounce paths can be found from Propositions 9
and 11 and the following two trigonometric identities: 
\begin{align}
\frac{\cos(\psi-\beta)}
    {\sin[\psi+(m-1)\beta]}
-\frac{\cos\psi}{\sin(\psi-m\beta)}
&   
+\sum_{k=1}^{2m-1}\frac{\cos(m\beta)\sin\beta}
    {\sin[\psi-(m-k)\beta]\sin[\psi-(m-k+1)\beta]} \nonumber \\\nonumber
&  \\\nonumber
& =2\sin(m\beta),
\end{align}
\[
\frac{\cos\psi}
    {\cos[(m-1)\beta]}  
-\sum_{k=1}^{m-1}\frac{\sin[\psi+(m-1)\beta]\sin\beta}
    {\cos[(m-k)\beta]\cos[(m-k-1)\beta]} =2\cos[\psi+(m-1)\beta].
\]

\vspace{0.2in}

\section{The Casimir Energy With Finite Top Plate}

Consider the case where the top plate is finite, lying
between $R=R_0$ and $R=R_1$ where $R_1>R_0$. Let the width of the top
plate be $W$ and the bottom plate be infinite in all
directions. 

\vspace{0.2in}

\setlength{\epsfxsize}{3in}
\begin{center}
   \mbox{\epsfbox{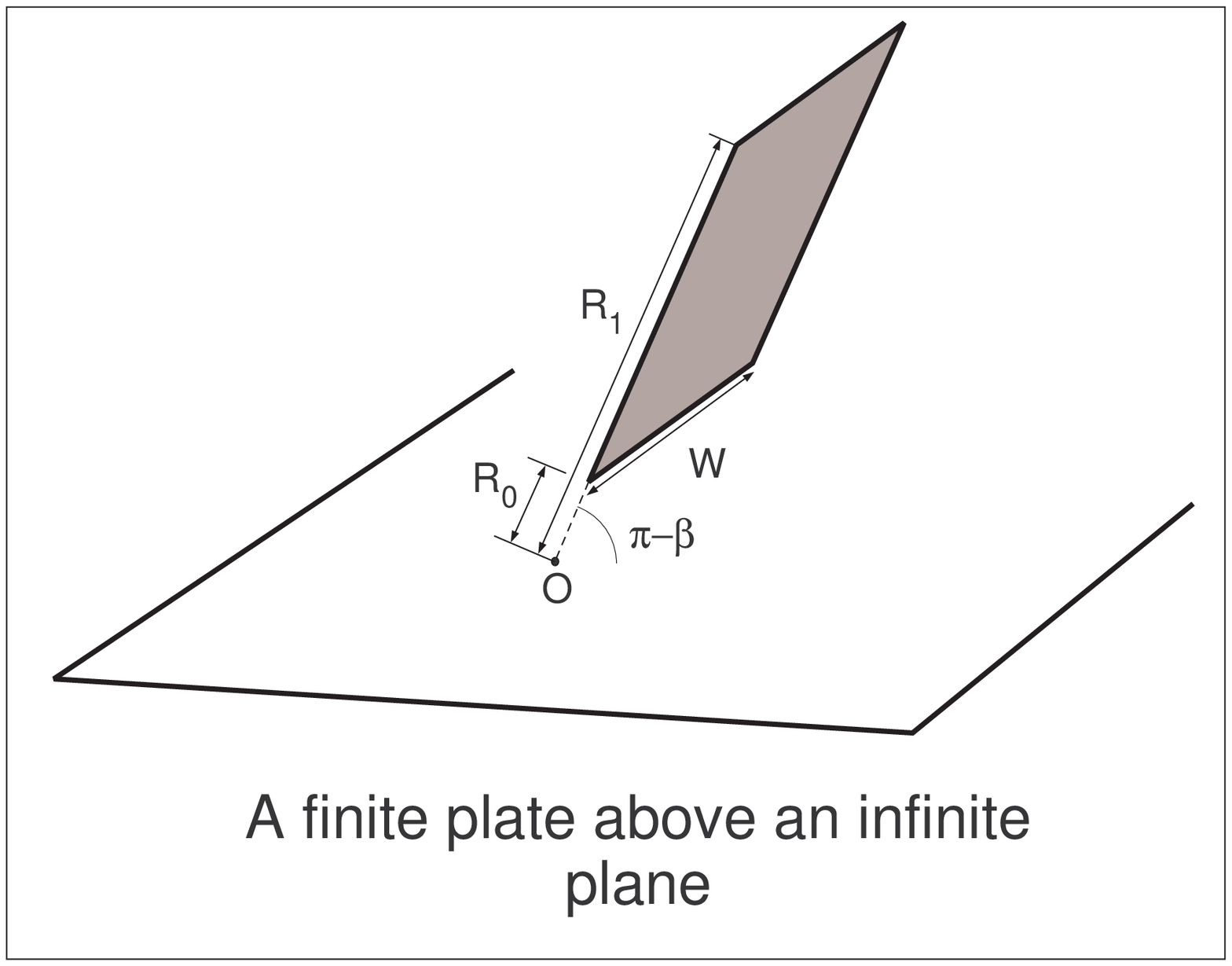}}
\end{center}

\vspace{0.2in}

The effect of limiting the upper plate is that some of the closed
paths must now be excluded from the energy as they run off the edge of
the plate. This restricts the domain of integration ${\cal{D}}_n$. 

Consider a closed $2m$- or $2m+1$-bounce path that first strikes the
horizontal plate. The closest and furthest points from the vertex have
$x^2$-coordinate:  

\vspace{0.2in}
\begin{center}
\begin{tabular}{r||c|c}\label{tab:rest1}
  & Closest & Furthest\\\hline\hline
m odd & $a_{m+1}$ & $a_{2m}$ \\\hline
m even & $a_{m}$ & $a_{2m}$ 
\end{tabular}
\end{center}
\vspace{0.2in}

For a $2m+1$-bounce path that strikes the top plate first, the closest
and furthest points from the vertex have $x^2$-coordinate: 

\vspace{0.2in}
\begin{center}
\begin{tabular}{r||c|c}\label{tab:rest2}
  & Closest & Furthest\\\hline\hline
m odd & $a_{m}$ & $a_{1}$ \\\hline
m even & $a_{m+1}$ & $a_{1}$ 
\end{tabular}
\end{center}
\vspace{0.2in}

We will use these to restrict the domains of integration, taking the
odd and even bounce cases separately.

\subsection{The Even Bounce Contribution}

We start by considering the case of a $2m$-bounce. 
The restrictions we must introduce for the finite plate, obtained
from (\ref{e:aevencl}) are

\vspace{0.2in}
\begin{center}
\begin{tabular}{r||c|c}
 m odd & $a_{m+1}\ge-R_0\cos\beta$ & $a_{2m}\le-R_1\cos\beta$ \\\hline
m even & $a_{m}\ge-R_0\cos\beta$ & $a_{2m}\le-R_1\cos\beta$ 
\end{tabular}
\end{center}
\vspace{0.2in}
Suppose that $m=2n$, then the restrictions are
\begin{equation}\label{e:circ1}
a_m=a_{2n}=\frac{R\cos(2n\beta)\cos\beta}{\sin(\psi-\beta)}\ge-R_0\cos\beta,
\end{equation}
\begin{equation}\label{e:circ2}
a_{2m}=a_{4n}=\frac{R\cos(2n\beta)\cos\beta}{\sin[\psi+(2n-1)\beta]}\le-R_1\cos\beta,
\end{equation}
which we put together to find
\[
-\frac{R_0\sin(\psi-\beta)}{\cos(2n\beta)}\le R\le-\frac{R_1\sin[\psi+(2n-1)\beta]}{\cos(2n\beta)}.
\]
We now consider the restriction this inequality places on $\psi$,
namely:
\[
-R_0\sin(\psi-\beta)\le-R_1\sin[\psi+(2n-1)\beta].
\]
Since $\beta-\frac{\pi}{2}\le\psi\le\frac{\pi}{2}$ and
$\beta>\frac{4n-1}{4n}\pi$ the function
\[
f(\psi)=-R_0\sin(\psi-\beta)+R_1\sin[\psi+(2n-1)\beta],
\]
is a decreasing function of $\psi$. Thus the inequality holds iff
$f(\frac{\pi}{2})\le0$, or  $-R_0\cos\beta\le -R_1\cos[(2n-1)\beta]$.

We consider two cases: either $f(\beta-\frac{\pi}{2})\le0$,
i.e. $R_0\le R_1\cos(2n\beta)$, or there exists
$\psi_0\in[\beta-\frac{\pi}{2},\frac{\pi}{2}]$ such that $f(\psi_0)=0$:
\[
-R_0\sin(\psi_0-\beta)=-R_1\sin[\psi_0+(2n-1)\beta].
\]
This can also be written
\begin{equation}\label{e:psi0even}
\tan\psi_0=\frac{R_0\sin\beta+R_1\sin[(2n-1)\beta]}
     {R_0\cos\beta-R_1\cos[(2n-1)\beta]}.
\end{equation}

In the former case the region of integration is
$[\beta-\frac{\pi}{2},\frac{\pi}{2}]$ while in 
the latter case it is $[\psi_0,\frac{\pi}{2}]$. 
The diagram below shows the regions of integration for a closed 4-bounce
for various wedge angles when $R_0=1$ and $R_1=2$. Here we can see the
transition from the latter to the former integration regions occuring
at 30$^0$ (note that, for illustration purposes, the scales on the two
axes are not equal).  

\vspace{0.2in}

\setlength{\epsfxsize}{4in}
\begin{center}
   \mbox{\epsfbox{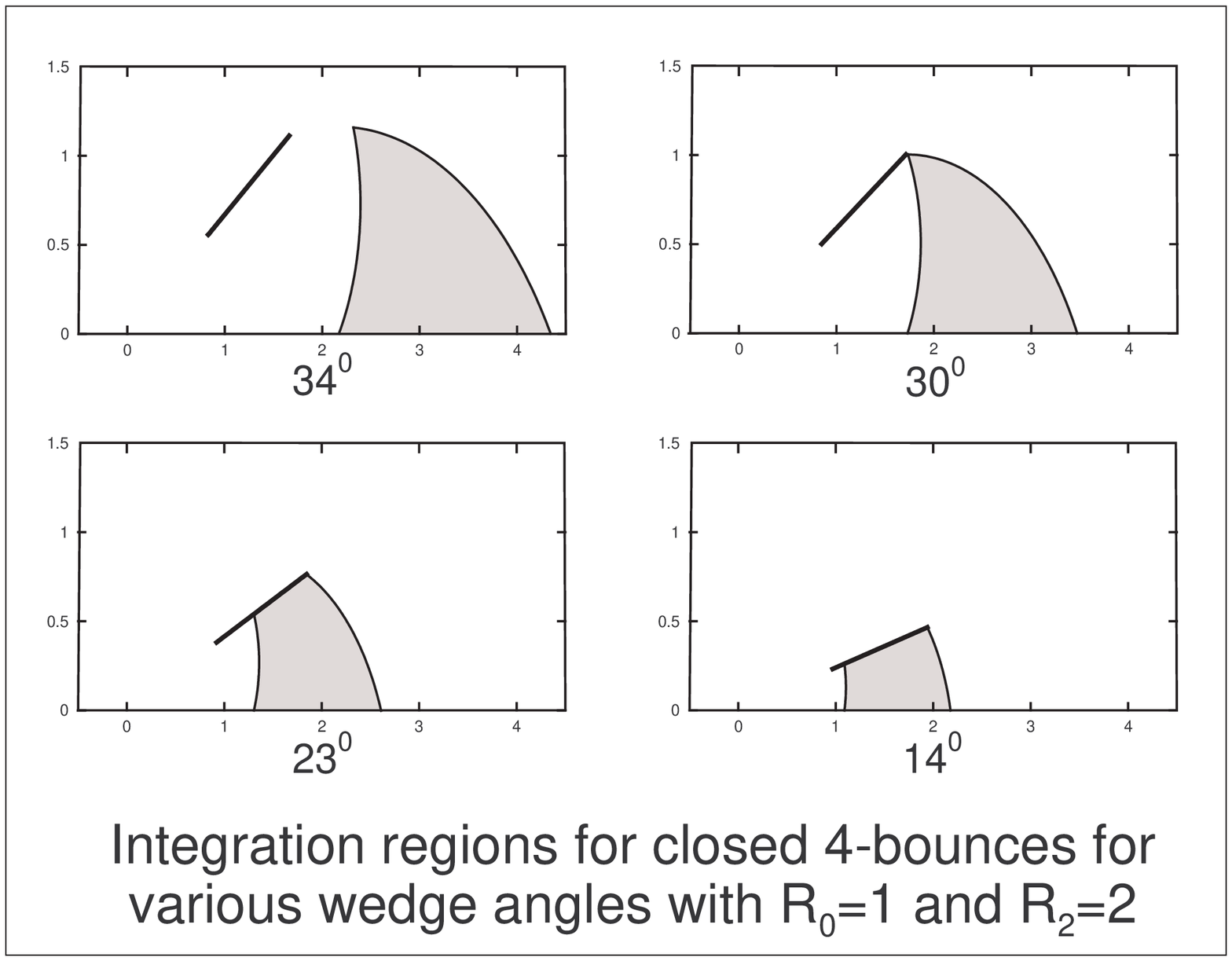}}
\end{center}

\vspace{0.2in}

In fact, it is easy to see that the regions of integration lie
between two circles which pass through the origin and whose second
point of intersection lies in the first quadrant at an angle $\psi_0$.
This follows from studying equations (\ref{e:circ1}) and (\ref{e:circ2}).
The region of integration is also seen to lie within the wedge.

If we let $\psi_1=\mbox{max}[\psi_0,\beta-\frac{\pi}{2}]$ then the
energy associated with a closed $4n$-bounce can be computed:

\begin{align}
{\cal{E}}_{4n}&=-\frac{\hbar cW}{32\pi^2}\int^{\frac{\pi}{2}}_{\psi_1}
   \int^{-\frac{R_1\sin(\psi+(2n-1)\beta)}{\cos(2n\beta)}}_{-\frac{R_0\sin(\psi-\beta)}{\cos(2n\beta)}}
\frac{1}{R^3\sin^4(2n\beta)}dR\;d\psi\nonumber\\\nonumber
&=\frac{\hbar cW}{64\pi^2\sin^4(2n\beta)}\int^{\frac{\pi}{2}}_{\psi_1}
\frac{\cos^2(2n\beta)}{R_1^2\sin^2[\psi+(2n-1)\beta]}-\frac{\cos^2(2n\beta)}{R_0^2\sin^2(\psi-\beta)}\;d\psi\\\nonumber
&=-\frac{\hbar cW\cos^2(2n\beta)}{64\pi^2\sin^4(2n\beta)}\left[
\frac{\cot[\psi+(2n-1)\beta]}{R_1^2}-\frac{\cot(\psi-\beta)}{R_0^2}
\right]^{\frac{\pi}{2}}_{\psi_1}\\\nonumber
&=-\frac{\hbar cW\cos^2(2n\beta)\cos\psi_1}{64\pi^2\sin^4(2n\beta)}\left(\frac{1}{R_0^2\cos\beta\sin(\psi_1-\beta)}-\frac{1}{R_1^2\cos[(2n-1)\beta]\sin[\psi_1+(2n-1)\beta]}\right)\nonumber.
\end{align}
Now suppose that $m=2n-1$. The restrictions 
together with (\ref{e:aevencl}), lead to
\[
a_{m+1}=a_{2n}=\frac{R\cos[(2n-1)\beta]\cos\beta}{\sin\psi}\ge-R_0\cos\beta,
\]
\[
a_{2m}=a_{4n-2}=\frac{R\cos[(2n-1)\beta]\cos\beta}{\sin[\psi+(2n-2)\beta]}\le-R_1\cos\beta,
\]
which we put together to find
\[
-\frac{R_0\sin\psi}{\cos[(2n-1)\beta]}\le R\le-\frac{R_1\sin[\psi+(2n-2)\beta]}{\cos[(2n-1)\beta]}.
\]
Thus this places the following inequality on $\psi$
\[
R_0\sin\psi\le R_1\sin[\psi+(2n-2)\beta],
\]
or
\begin{equation}\label{e:psi0odd}
\tan\psi_0=\frac{R_1\sin[(2n-2)\beta]}
     {R_0-R_1\cos[(2n-2)\beta]}.
\end{equation}

A similar monotone argument to the $m=2n$ case shows that the
inequality holds iff $R_0\le R_1\cos[(2n-2)\beta]$.

We consider two cases: either  
\[
-R_0\cos\beta\le -R_1\cos[(2n-1)\beta],
\]
or there exists $\psi_0\in[\beta-\frac{\pi}{2},\frac{\pi}{2}]$ such
that
\[
R_0\sin\psi_0=R_1\sin[\psi_0+(2n-2)\beta].
\]
In the former case the region of integration is
$[\beta-\frac{\pi}{2},\frac{\pi}{2}]$ while in 
the latter case it is $[\psi_0,\frac{\pi}{2}]$. As before, these lie
between two circles which intersect at the origin and at a
point in the first quadrant. 

Letting
$\psi_1=\mbox{max}[\psi_0,\beta-\frac{\pi}{2}]$ the resulting energy
integrates up to
\[
{\cal{E}}_{4n-2}=-\frac{\hbar cW\cos^2[(2n-1)\beta]\cos\psi_1}{64\pi^2\sin^4[(2n-1)\beta]}\left(
\frac{1}{R_0^2\sin\psi_1}-\frac{1}{R_1^2\cos[(2n-2)\beta]\sin[\psi_1+(2n-2)\beta]}\right).
\]

We combine these results for the even bounce case:

\begin{Prop}\label{p:fineven}
Given $\beta$, $R_0$ and $R_1$, define $m_0$ and $m_1$ by 
$\frac{2m_0-1}{2m_0}\pi<\beta\le\frac{2m_0+1}{2m_0+2}\pi$, and either 
$\cos(m_1\beta)\le\frac{R_0}{R_1}\le\frac{\cos[(m_1-1)\beta]}{\cos\beta}$ 
for $m_1$ even or
$\frac{\cos(m_1\beta)}{\cos\beta}\le\frac{R_0}{R_1}\le\cos[(m_1-1)\beta]$ 
for $m_1$ odd. 

Then the even contribution to the Casimir energy is
\[
{\cal{E}}_{\mbox{even}}=2\sum_{m=1}^{m_0}{\cal{E}}_{2m}^1, 
\]
when $m_0\le m_1$ and
\[
{\cal{E}}_{\mbox{even}}=2\sum_{m=1}^{m_1-1}{\cal{E}}_{2m}^1+2{\cal{E}}_{2m_1}^0 
\]
when $m_0> m_1$, where
\begin{equation}\label{e:even1}
{\cal{E}}_{2m}^1=\frac{\hbar cW\cos^2(m\beta)\sin\beta}{64\pi^2\sin^4(m\beta)}\left(\frac{1}{R_0^2\cos\beta}-\frac{1}{R_1^2\cos[(m-1)\beta]\cos(m\beta)}\right),
\end{equation}
and
\begin{equation}\label{e:even0even}
{\cal{E}}_{2m}^0=-\frac{\hbar
  cW\cos^2(m\beta)(R_0\cos\beta-R_1\cos[(m-1)\beta])^2}{64\pi^2\sin^4(m\beta)\cos\beta\cos[(m-1)\beta]\sin(m\beta)R_0^2R_1^2},  
\end{equation}
for $m_1$ even, and
\begin{equation}\label{e:even0odd}
{\cal{E}}_{2m}^0=-\frac{\hbar
  cW\cos^2(m\beta)(R_0-R_1\cos[(m-1)\beta])^2}{64\pi^2\sin^4(m\beta)\sin[(m-1)\beta]\cos[(m-1)\beta]R_0^2R_1^2},
\end{equation}
for $m_1$ odd.
\end{Prop}

\vspace{0.2in}

\noindent{\bf Note}:
\begin{enumerate}
\item[(1)] If the above restrictions hold for $n=n_0$ then they hold
  for $n=n_0-1$. In other words, the restrictions give an upper bound
  on the number of contributions of closed bounces.
\item[(2)] As $\beta\rightarrow\pi$ the restrictions reduce to
  $\frac{2m_0-1}{2m_0}\pi<\beta\le\frac{2m_0+1}{2m_0+2}\pi$. 
\end{enumerate}

\vspace{0.2in}

\subsection{The Odd Bounce Contribution}

The situation for odd bounce paths is quite different than for even
bounce paths. Due to the finite size of the top plates there are again
restrictions on the regions of integration. However, rather than
eliminating the contribution from higher bounces, as in the even case, these contributions 
give lower bounce paths originating from regions outside of the wedge.

In the diagram below the region ${\cal{D}}_5$ for the closed 5-bounce
paths is shown for $\gamma=22.5^0$, $R_0=1$ and $R_1=2$. This
includes a region outside of the wedge (darker shading), which is the
reflection of the 7-bounce region that would exist if the top plate
were larger.  

\vspace{0.2in}

\setlength{\epsfxsize}{4in}
\begin{center}
   \mbox{\epsfbox{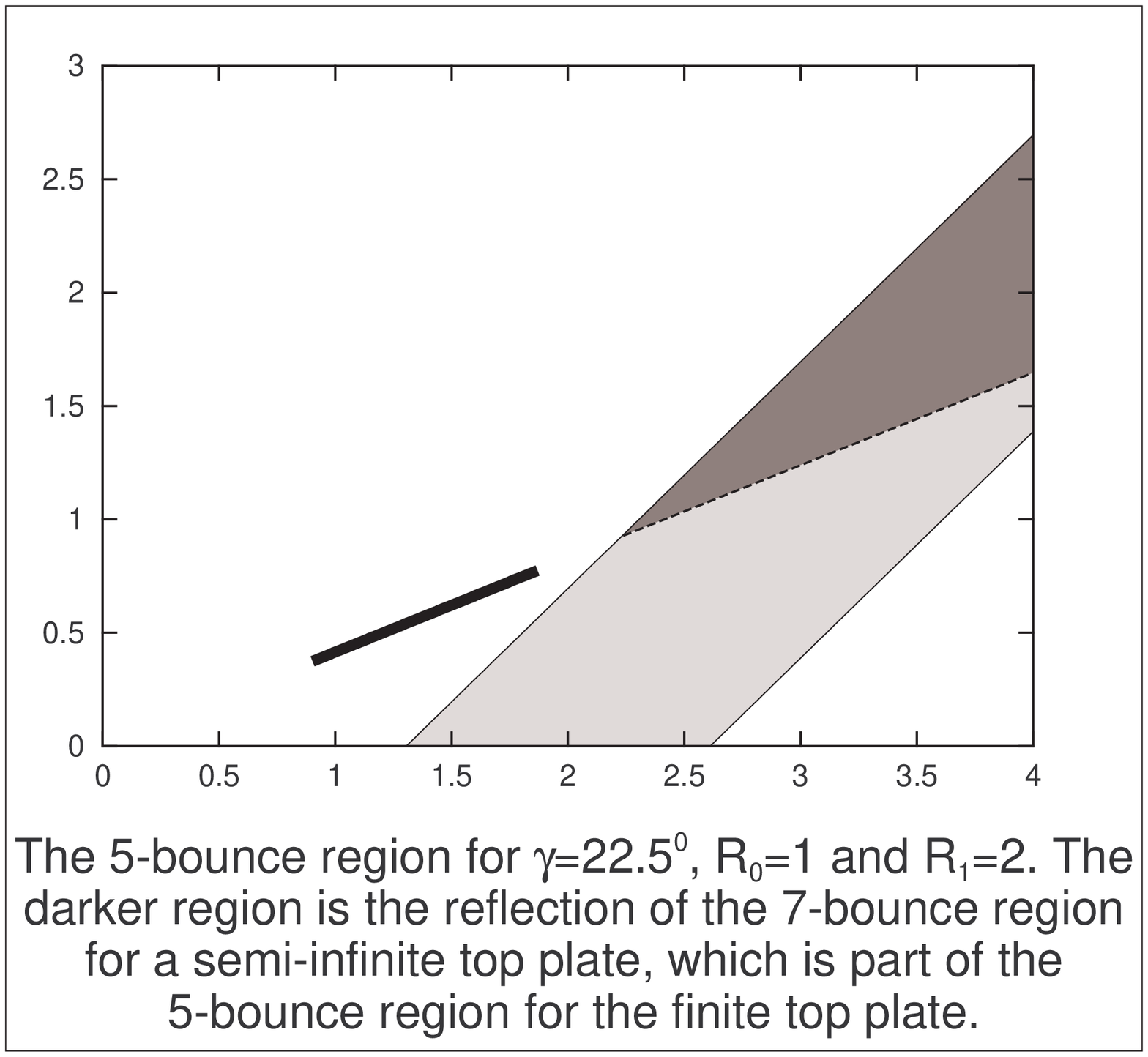}}
\end{center}

\vspace{0.2in}

Since reflection preserves lengths and areas, the total contribution
of the odd bounces can be calculated by integrating with only the
lower bound on $R$. That is, we compute the energy for a semi-infinite top plate.

Consider a closed $2m+1$-bounce path, with $m=2n$, which
hits the horizontal plate first. As we saw, such a path will exist
only if $\beta>\frac{2m-1}{2m}\pi$.  The sequence of intersection points
(\ref{e:aoddcl}), together with the
restriction leads to the lower bound on $R$:
\[
R\ge-\frac{R_0\cos\beta}{\sin(\psi+2n\beta)}.
\]

On the other hand, we have seen that a closed $4n+1$ bounce which hits
the horizontal plate first exists
only if $\psi>2n(\pi-\beta)$. Thus the limits of integration for
$\psi$ are $\psi_0\le\psi\le\frac{\pi}{2}$, where
$\psi_0=\mbox{max}[2n(\pi-\beta),\beta-\frac{\pi}{2}]$. Thus

\begin{align}
{\cal{E}}_{4n+1}^{H}&=\frac{\hbar cW}{32\pi^2}\int^{\frac{\pi}{2}}_{\psi_0}\int^{\infty}_{-\frac{R_0\cos\beta}{\sin(\psi+2n\beta)}}
\frac{1}{R^3\cos^4(\psi+2n\beta)}dR\;d\psi\nonumber\\\nonumber
&=\frac{\hbar cW}{64\pi^2\cos^2\beta R_0^2}\int^{\frac{\pi}{2}}_{\psi_0}\frac{\sin^2(\psi+2n\beta)}{\cos^4(\psi+2n\beta)}\;d\psi\\\nonumber
&=-\frac{\hbar cW}{192\pi^2\cos^2\beta R_0^2}\left[\frac{\sin^3(\psi+2n\beta)}{\cos^3(\psi+2n\beta)}\right]^{\frac{\pi}{2}}_{\psi_0}\;d\psi\\\nonumber
&=-\frac{\hbar cW}{192\pi^2\cos^2\beta R_0^2}\left(\cot^3(2n\beta)-\epsilon\cot^3[(2n+1)\beta]\right),\nonumber
\end{align}
where $\epsilon=0$ if $2n(\pi-\beta)>\beta-\frac{\pi}{2}$ and
$\epsilon=1$ if $2n(\pi-\beta)\le\beta-\frac{\pi}{2}$. We denote these
two contributions by ${\cal{E}}_{4n+1}^{H\epsilon}$.

Consider the case where $m=2n-1$. Again we encounter a restriction on
the domain of integration analogous to the previous case. The result
is:
\[
{\cal{E}}_{4n-1}^{H\epsilon}=-\frac{\hbar cW}{192\pi^2R_0^2}
   \left(\cot^3[(2n-1)\beta]-\epsilon\cot^3(2n\beta)
   \right),
\]
where $\epsilon=0$ when $(2n-1)(\pi-\beta)>\beta-\frac{\pi}{2}$ and
$\epsilon=1$ when $(2n-1)(\pi-\beta)<\beta-\frac{\pi}{2}$. 

Finally we look at odd bounces that first strike the top plate. The
sequence of intersection points is now given by (\ref{e:aoddclt}). The
resulting energies are: 

\[
{\cal{E}}_{4n+1}^{T\epsilon}=-\frac{\hbar cW}{192\pi^2R_0^2}
   \left(\cot^3(2n\beta)-\epsilon\cot^3[(2n+1)\beta]\right),
\]
\[
{\cal{E}}_{4n-1}^{T\epsilon}=-\frac{\hbar
  cW}{192\pi^2\cos^2\beta  R_0^2}\left(\cot^3[(2n-1)\beta]-\epsilon\cot^3(2n\beta)\right), 
\]
where $\epsilon=0$ when $m\beta-(2m-1)\pi>\frac{\pi}{2}$ and $\epsilon=1$
when $m\beta-(2m-1)\pi<\frac{\pi}{2}$. 

Summing these four energies we have that the contribution of the odd
bounces to the Casimir energy in a wedge with
$\frac{2m_0-1}{2m_0}\pi<\beta\le\frac{2m_0+1}{2m_0+2}\pi$ is
\[
{\cal{E}_{\mbox{odd}}}=\sum_{m=1}^{m_0-1}\left({\cal{E}}_{2m+1}^{H1}+{\cal{E}}_{2m+1}^{T1}\right)+{\cal{E}}_{2m_0+1}^{H0}+{\cal{E}}_{2m_0+1}^{T0}=
-\frac{\hbar cW(1+\cos^2\beta)\cos\beta}{192\pi^2R_0^2\sin^3\beta}
\]

Note that, since this term does not involve $R_1$, it is independent of
the length of the top plate. From a physical point of view, such
contributions to the force have no significance, and are excluded
\cite{jas}. Thus we only have the even bounce contributions, and Main
Theorem 2 follows from Proposition \ref{p:fineven} with
$\gamma=\pi-\beta$.

\section{The Parallel Plate Limit}

From Casimir's original work \cite{cas}, the Dirichlet energy between
two parallel plates of area $A$ and separation $L$ was computed to be
\[
{\cal{E}}=-\frac{\hbar c\pi^2A}{1440\;L^3}.
\]
To date, this is the only closed analytic expression for the Casimir
energy. In this section we retrieve this result as the limit of our
expressions for the energy between non-parallel plates. 

Let us consider the limit $\beta\rightarrow\pi$. As it stands, the
Casimir energy between a finite plate and an infinite plane, as given
above, diverges as $\beta\rightarrow\pi$. This is because the separation
between the boundaries goes to zero in this limit. 

Before taking the limit, we fix the non-zero minimum separation $L$
between the plates by letting $R_0\sin\beta=L$ and
$b=R_1-R_0$, where $b$ is the length of the top plate (see the diagram
below).

\vspace{0.2in}

\setlength{\epsfxsize}{4in}
\begin{center}
   \mbox{\epsfbox{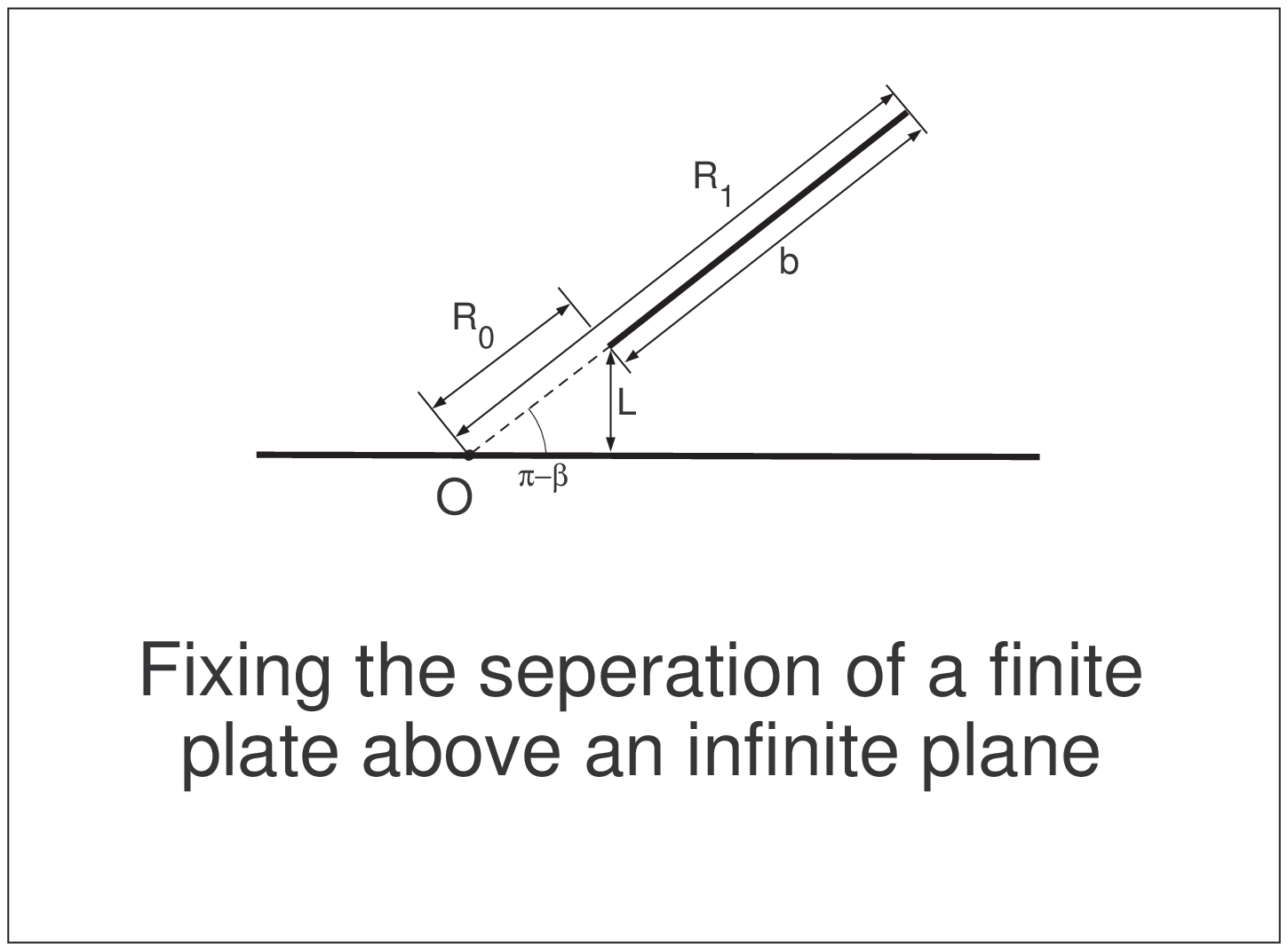}}
\end{center}

\vspace{0.2in}

While for each $\beta<\pi$ there is only a finite number of
contributions to the Casimir energy, as $\beta\rightarrow\pi$ we pick
up an infinite number of terms in the sum. Thus the parallel plate limit
needs careful consideration. 

\vspace{0.2in}

\noindent {\bf Main Theorem 3}.

{\it
In the parallel plate limit we retrieve Casimir's orginal result:
\[
\lim_{\beta\rightarrow\pi}{\cal{E}}=-\frac{\hbar c\pi^2bW}{1440\;L^3}.
\]
}

\begin{pf}
First we interchange the limit with the sum:
\[
\lim_{\beta\rightarrow\pi}{\cal{E}}=2\sum_{m=1}^{\infty}\lim_{\beta\rightarrow\pi}{\cal{E}}_{2m}.
\]

We can do this since for each angle
$\beta<\pi$ there are only a finite number of non-zero contributions, and
the energy is therefore equicontinuous.

Next, the limit can be computed directly.
As noted earlier, the extra restrictions on $\psi$ 
encountered in computing the even contributions for the finite plate
disapppear as $\beta\rightarrow\pi$. Thus, in the limit,
\[
{\cal{E}}_{2m}\rightarrow\frac{\hbar cW\cos^2(m\beta)\sin^3\beta}{64\pi^2\sin^4(m\beta)}\left(\frac{1}{L^2\cos\beta}-\frac{1}{(L+b\sin\beta)^2\cos[(m-1)\beta]\cos(m\beta)}\right).
\]
We apply L'H\^opital's four times to show that 
\[
\lim_{\beta\rightarrow\pi}\frac{\sin^3\beta((L+b\sin\beta)^2\cos[(m-1)\beta]\cos(m\beta)-L^2\cos\beta)}{\sin^4(2m\beta)}=-\frac{2bL}{m^4}.
\]
Thus
\[
\lim_{\beta\rightarrow\pi}{\cal{E}}_{2m}=-\frac{\hbar cbW}{32\pi^2L^3\;m^4}
\]

Finally the result follows:
\[
\lim_{\beta\rightarrow\pi}{\cal{E}}=2\sum_{m=1}^{\infty}\lim_{\beta\rightarrow\pi}{\cal{E}}_{2m}=-\frac{\hbar
  cbW}{16\pi^2L^3}\sum_{m=1}^{\infty}\frac{1}{m^4}=-\frac{\hbar c\pi^2bW}{1440\;L^3}.
\]

\end{pf}

\end{document}